\documentclass[a4paper,11pt,twoside]{article}
\usepackage{pst-fill,pst-grad}  
\usepackage{textcomp}  
\usepackage[english]{babel}  
\usepackage[utf8]{inputenc} 
\usepackage{graphicx} 
\usepackage{amsmath}
\usepackage{float} 
\usepackage{fancyhdr}
\usepackage[matrix,arrow,curve]{xy}
\usepackage{pstricks} 
\usepackage{amsmath,amsfonts,verbatim,afterpage,theorem,euscript,mathrsfs,amssymb} 
\usepackage{amsfonts}
\usepackage{amssymb} 
\usepackage{array}   
\usepackage{dsfont} 
\usepackage[titletoc]{appendix}
\usepackage[colorlinks=true,linkcolor=blue,citecolor=red]{hyperref}
\usepackage{authblk} 
\usepackage{color}
 
\newtheorem{Definition}{Definition}[section] 
\newtheorem{Proposition}{Proposition}[section]

\newtheorem{Lemme}{Lemma}[section]
\newtheorem{Theoreme}{Theorem}[section]

\newtheorem{Corollaire}{Corollary}[section]
\newtheorem{Remarque}{Remark}[section]
\def \f{\vec{f}} 
\def \g{\vec{g}}

\def \vu{\vec{u}}
\def \P{\mathbb{P}}
\def \U{\vec{U}}

\def \V{\vec{V}}
\def \v{\vec{v}}

\def \B{\vec{B}}
\def \b{\vec{b}}
\def \bV{\textbf{V}}
\def \Rt{\mathbb{R}^{3}}

\def \R{\mathbb{R}}

\def \finpv{\hfill $\blacksquare$}
\def \pv{{\bf{Proof.}}~}

\def \ds{\displaystyle}
 
\def \bl{\textcolor{black}} 
\title{ \bf Some remarks on the  regularity of  weak  solutions for the stationary Ericksen-Leslie and  MHD systems}  

\author[1]{ Oscar Jarr\'in \footnote{corresponding author: oscar.jarrin@udla.edu.ec}}  
\affil[1]{\scriptsize Escuela de Ciencias Físicas y Matemáticas, Universidad de Las Américas, Vía a Nayón, C.P.170124, Quito, Ecuador.}  

\begin{document} 

\maketitle  

\begin{abstract}
\textbf{Abstract.} We consider  two  elliptic coupled systems  of relevance in   the fluid dynamics. These systems   are posed on the whole space $\Rt$ and they consider the action of external forces. The first system deals with the simplified Ericksen-Leslie (SEL) system, which describes the dynamics of liquid crystal flows. The second system  is  the time-independent magneto-hydrodynamic (MHD) equations. For the (SEL) system,  we obtain  a new criterion to  improve the regularity of weak solutions, provided that they belong to  some  homogeneous Morrey space. As a bi-product, we also obtain some new regularity criterion for the stationary Navier-Stokes equations and for a nonlinear harmonic map flow.  This new regularity criterion  also holds true for the (MHD) equations. Furthermore,  for this last  system  we are able to use the Gevrey class  to prove that all finite energy weak solutions are analytic functions, provided the external forces belong to some Gevrey class. \\[5mm]
[3mm]
\end{abstract}
\section{Introduction} 
This  note  considers with two elliptic  coupled systems  in  the fluid dynamics.  The first system arises from the  study of the dynamics in liquid crystal flows. This system is posed on the whole space $\Rt$ and strongly couples the  incompressible  and stationary (time-independent) Navier-Stokes equations  with a nonlinear harmonic map flow as follows: 
\begin{equation}\label{EickLes}  
\left\{ \begin{array}{ll}\vspace{2mm} 
-\Delta \U + \text{div} (\U \otimes \U)   + \text{div} ( \vec{\nabla} \otimes  \V \odot \vec{\nabla} \otimes \V) + \vec{\nabla} P= \text{div}(\mathbb{F}), \\ \vspace{2mm} 
-\Delta \V + \text{div}(\V \otimes \U)   - \vert  \vec{\nabla} \otimes \V  \vert^2\, \V = \text{div}(\mathbb{G}). \\
\text{div}(\U)=0. \end{array} \right. 
\end{equation}
Here, $\ds{\vec{\nabla} \otimes  \V =(\partial_i V_j)_{1\leq i,j\leq 3}}$, denotes the deformation tensor of the vector field $\V$ and moreover, for $i=1,2,3$, each component of the vector field $\ds{ \text{div} ( \vec{\nabla} \otimes  \V \odot \vec{\nabla} \otimes \V) }$ explicitly writes down as:  
\begin{equation}\label{super-critical}
\left[\text{div} ( \vec{\nabla} \otimes  \V \odot \vec{\nabla} \otimes \V)\right]_i= \sum_{j=1}^{3} \sum_{k=1}^{3} \partial_{j}( \partial_i V_k \, \partial_j V_k).
\end{equation} 
The velocity of the fluid $\U:\Rt \to \Rt$ and the pressure $P:  \Rt \to \mathbb{R}$  are the classical unknowns of the fluid mechanics. Moreover, this system also considers a third unknown $\V  :\Rt  \to \mathbb{S}^{2}$, where $\mathbb{S}^{2}$ denotes the unit sphere in $\Rt$. The \emph{unit vector field} $\V$  represents the macroscopic orientation of the nematic liquid crystal molecules \cite{Les}.  We also take into account the action of external forces,  which according to \cite{Gennes}  can be written  as the divergence of the tensors $\mathbb{F}=(F_{i,j})_{1\leq i,j \leq 3}$ and $\mathbb{G}=(G_{i,j})_{1\leq i,j\leq 3}$, with $F_{i,j}, G_{i,j}: \Rt \to \R$. Finally, the equation $\text{div}(\U)=0$ always represents the fluid's incompressibility.

\medskip

The elliptic  system (\ref{EickLes}) is the time-independent  counterpart  of the following parabolic (time-dependent) system:   
\begin{equation}\label{EickLes-non-stat}
\left\{ \begin{array}{ll}\vspace{3mm} 
\partial_t \vu -\Delta \vu + \text{div}(\vu \otimes \vu)  + \text{div} ( \nabla \otimes  \v \odot  \nabla  \otimes  \v) + \nabla p=\text{div}(\mathbb{F}), \quad \text{div}(\vu)=0, \\ 
\partial_t \v-\Delta \v + \text{div}(\v \otimes \vu)  - \vert  \nabla \otimes \v  \vert^2 \, \v =\text{div}(\mathbb{G}),    
\end{array} \right. 
\end{equation} 
also known as the \emph{ simplified Ericksen-Leslie system}. This parabolic system was proposed by H.F. Lin in \bl{\cite{Lin}} as a simplification of the general \emph{Ericksen-Leslie system} which models the hydrodynamic flow of nematic liquid crystal material \cite{Ericksen,Les}. The simplified Ericksen-Leslie system  has been successful to model the dynamical behavior of nematic liquid crystals. More precisely, it provides a good macroscopic description  of the evolution of the material under the influence of the fluid velocity field,  and moreover, it provides a good   macroscopic description  of rod-like liquid crystals. See the book \bl{\cite{Gennes}} for more details. 

\medskip

The  system (\ref{EickLes-non-stat}) has recently attired the interest in the research community of mathematical fluid dynamics. It is worth mentioning  one of the major challenges is due, on the one hand,  to the  presence of the trilinear term $\vert  \vec{\nabla} \otimes \v  \vert^2 \, \v $ in the second equations of this system and, on the other hand, to the presence of  the super-critical nonlinear term  $ \text{div} ( \nabla \otimes  \v \odot  \nabla  \otimes  \v)$ defined in (\ref{super-critical}). Precisely,  the double derivatives in this last  term make it  more  delicate to treat than the classical nonlinear transport term: $\ds{\text{div}(\vu \otimes \vu)}$. These  facts make challenging  the study of both (\ref{EickLes})  and (\ref{EickLes-non-stat}). See, for instance, the articles \cite{HaLiZ,OJ,LinWang,LinLiu,LinLinWang,Wang} and the references therein.

\medskip

{\bf Some previous works in the homogeneous case}.  When $\mathbb{F}=0$ and $\mathbb{G}=0$, the first works on the mathematically  study of the system (\ref{EickLes-non-stat})   were devoted to the  existence of global in time weak solutions  \cite{LinWang,LinLinWang}. Thereafter,  in the spirit of the celebrated result by  H. Koch \& D. Tataru \cite{KochTataru}, the global well-posedness of small solutions in the space $BMO^{-1}(\Rt)$ was proven 
in \cite{Wang}. 

\medskip

Concerning some regularity issues, T. Huang proved in  \cite{Huang} a regularity criterion on weak solutions  of (\ref{EickLes-non-stat})  in the framework of the Lebesgue spaces.  This result also holds true for the stationary system (\ref{EickLes}) in the case $\mathbb{F}=\mathbb{G}=0$. Indeed, first  we consider  a weak solution of (\ref{EickLes}) as a couple $(\U, \V)$  where $\ds{(\U,  \vec{\nabla} \otimes \V)  \in H^1(\Rt)}$. Thereafter, we obtain $\U \in  \mathcal{C}^{\infty}(\Rt)$ and  $ \V \in  \mathcal{C}^{\infty}(\Rt)$, provided that $\U  \in  L^p(\Rt)$ and  $\vec{\nabla} \otimes  \V \in  L^p(\Rt)$, with $p>3$. 

\medskip 

Regularity of weak solutions is an important question to a better mathematically comprehension of the system (\ref{EickLes}), and moreover, it is also one of the key assumptions  when studying another relevant problem of this system in the homogeneous case: the uniqueness of weak solutions. Precisely, when $\mathbb{F}=\mathbb{G}=0$ we always have the trivial solution $(\U, P,\V)=(0,0,0)$ and  we  look for some functional spaces in which this solution is the unique one. This problem, also known the \emph{Liouville-type problem},  was recently studied in \cite{OJ} where the main interest is the use of more general spaces than the $L^p$ spaces, for instance, the Lorentz and the Morrey spaces. However, to the best of our knowledge, the regularity of weak solutions in these spaces was not studied before and it must be assumed.  

\medskip 

{\bf A new regularity criterion in the non-homogeneous case.}  Motivated by this last question,  in this note  we study some  new \emph{a priori} conditions in the setting of the Morrey spaces  to improve the regularity of weak solutions of the system (\ref{EickLes}). Moreover, we also consider the more general case under the action of the external forces $\text{div}(\mathbb{F})$ and $\text{div}(\mathbb{G})$.  

\medskip

We shall consider here a fairly general notion of weak solutions, which is given in the following: 
\begin{Definition}\label{Def-weakSol-EL} Let $\mathbb{F}, \mathbb{G} \in \mathcal{D}'(\Rt)$.  A  weak solution of the coupled system (\ref{EickLes}) is  a triplet $(\U, P, \V)$,  where:  $\U \in L^{2}_{loc}(\Rt)$, $P \in \mathcal{D}'(\Rt)$, $\vert \V(x)\vert =1$ for almost all  $x\in \Rt$  and $\vec{\nabla} \otimes \V \in L^{2}_{loc}(\Rt)$, such that  it verifies (\ref{EickLes})  in the distributional sense. 
\end{Definition}	
It is worth observing we use minimal conditions on the functions $\U, \V$ and $P$  to ensure that  all the terms in  (\ref{EickLes})  are well defined as distributions. Moreover, we let the pressure $P$ to be a very general object as we only have $P\in \mathcal{D}'(\Rt)$.  

\medskip 

As $\U$ and $\vec{\nabla} \otimes \V$ are locally  square integrable functions, in order to improve their regularity  we look for some natural conditions on the local quantities  $\ds{\int_{B(x_0,R)} \vert \U(x) \vert^2\, dx }$ and $\ds{\int_{B(x_0,R)} \vert \vec{\nabla} \otimes \V(x) \vert^2\, dx }$,  where $B(x_0, R)$ denotes the ball of center  $x_0 \in \Rt$ and radius $R>0$. Thus,  the   Morrey spaces appear naturally. We recall that for  a parameter $2 < p <+\infty$, the homogeneous Morrey space $\dot{M}^{2,p}(\Rt)$  is the Banach space of functions $f \in L^{2}_{loc}(\Rt)$ such that:
\begin{equation}\label{Def-Morrey}
\Vert f \Vert_{\dot{M}^{2,p}}= \sup_{x_0 \in \Rt,\,  R>0}\,  R^{\frac{3}{p}} \left(  \frac{1}{\vert B(x_0, R) \vert} \int_{B(x_0,R)} \vert f (x) \vert^2 dx \right)^{\frac{1}{2}}  <+\infty,
\end{equation} 
where $\vert  B(x_0, R) \vert \simeq R^3$ is the Lebesgue measure of the ball.  The parameter $p$ measures  the decaying rate of the (local) mean quantity $\ds{\left(  \frac{1}{\vert B(x_0,R)\vert} \int_{B(x_0,R)} \vert f (x) \vert^2 dx \right)^{\frac{1}{2}}}$ as  $R$ goes to infinity.  This is a homogeneous space of degree $-\frac{3}{p}$, and moreover,  we have the following chain of continuous embeddings  $\ds{L^{p}(\Rt)\subset L^{p,q}(\Rt)  \subset \dot{M}^{2,p}(\Rt)}$.  Here,    $L^{p,q}(\Rt)$ (with $p<q\leq +\infty$) denotes a Lorentz space which  describes the decaying properties of functions in a different setting. See the book \bl{\cite{DCh}} for a detailed study of these spaces. 

\medskip

Finally, for the parameter $p>2$  given above, and for the regularity parameter $k\geq 0$, we introduce now the Sobolev-Morrey space  
\[ \mathcal{W}^{k,p}(\Rt) = \left\{ f \in \dot{M}^{2,p}(\Rt) : \,\, \partial^{\alpha} f \in  \dot{M}^{2,p}(\Rt), \,\, \text{for all multi-indice}\,\, \vert \alpha \vert \leq k   \right\}. \]
Moreover, we  denote by $W^{k,\infty}(\Rt)$ the classical Sobolev space of bounded functions with bounded weak derivatives until the order $k$.  Then, our first result reads as follows: 

\begin{Theoreme}\label{Th} Let  $\big( \U, P, \V \big)$ be a  weak solution of the coupled system (\ref{EickLes}) given in Definition \ref{Def-weakSol-EL}. We assume  $\U \in \dot{M}^{2,p}(\Rt)$ and $\vec{\nabla}\otimes \V \in \dot{M}^{2,p}(\Rt)$, with $p>3$. Then, if for  $k \geq 0$:  
	\begin{equation}\label{Cond-Th-1}
	\mathbb{F}, \, \mathbb{G} \in \mathcal{W}^{k+1,p}(\Rt) \cap W^{k+1,\infty}(\Rt),
	\end{equation}	
	it follows that  $\U \in \mathcal{W}^{k+2,p}(\Rt)$, $P \in \mathcal{W}^{k+1,p}(\Rt)$ and $\V \in \mathcal{W}^{k+2,p}(\Rt)$.  Moreover,   for all multi-indice $\vert \alpha  \vert \leq k+1$, the functions $\partial^{\alpha} \U$  and $\partial^{\alpha} \V$  are H\"older continuous with exponent $0<1-3/p<1$, while for $\vert \alpha \vert \leq k$ the function $\partial^{\alpha} P$  is also  H\"older continuous with the same exponent. 			
\end{Theoreme}	 

\begin{Remarque} Recall that the external forces acting on the system (\ref{EickLes})  are given by $\text{div}(\mathbb{F})$ and $\text{div}(\mathbb{G})$. Then,  by (\ref{Cond-Th-1}) we have $\text{div}(\mathbb{F}), \, \text{div}(\mathbb{G}) \in \mathcal{W}^{k,p}(\Rt)$ which yields a gain  of regularity of weak  solutions of the order $k+2$.  This (expected) maximum  gain of regularity is given by the effects of the Laplacian operator  in the system (\ref{EickLes}).  We refer to Remark \ref{Rmk-tech} for the technical details. 
\end{Remarque}	 
\begin{Remarque}\label{Rmk} For the particular  homogeneous case when $\mathbb{F}=\mathbb{G}=0$, we obtain that weak solutions  of the system (\ref{EickLes}) verify  $(\U, P, \V) \in \mathcal{C}^{\infty}(\Rt)$, provided that $\U \in \dot{M}^{2,p}(\Rt)$ and $\vec{\nabla}\otimes \V \in \dot{M}^{2,p}(\Rt)$, with $p>3$. As explained, this particular result is of interest in connection to the Liouville-type problem for (\ref{EickLes}) in the Morrey spaces \cite{OJ}. 
\end{Remarque}	

Mathematically, the coupled system (\ref{EickLes}) is also of interests as it contains  two relevant equations. On the one hand, by setting  $\V$ a constant unitary vector we  get the stationary  and forced Navier-Stokes equations:
\begin{equation}\label{NS} 
-\Delta \U + \text{div} (\U \otimes \U)  + \vec{\nabla} P= \text{div}(\mathbb{F}), \qquad \text{div}(\U)=0.
\end{equation}
On the other hand, by setting now $\U=0$ we obtain the following harmonic map flow: 
\begin{equation}\label{Harmonic-map}
-\Delta \V   - \vert  \vec{\nabla} \otimes \V  \vert^2\, \V = \text{div}(\mathbb{G}).
\end{equation}
As a direct consequence of this result we obtain a new regularity criterion for these equations:
\begin{Corollaire}
\begin{enumerate}
\item[] 	
\item[1.] Let $(\U, P) \in L^{2}_{loc}(\Rt)\times \mathcal{D}'(\Rt)$ be a  weak solution of the equation (\ref{NS}). If $\U \in \dot{M}^{2,p}(\Rt)$, with $p>3$, and $\mathbb{F}$ verifies (\ref{Cond-Th-1}) with $k\geq 0$, then we have $\U \in \mathcal{W}^{k+2,p}(\Rt)$ and $P \in \mathcal{W}^{k+1,p}(\Rt)$. 	 Furthermore, for  $\vert \alpha \vert \leq k+1$ the function $\partial^\alpha \U$ is $(1-3/p)-$H\"older continuous, while this holds true for $\partial^\alpha P$ with $\vert \alpha \vert \leq k$.

\medskip 

\item[2.] Let $\V(x) \in \mathbb{S}^2$ with $\vec{\nabla} \otimes \V \in L^{2}_{loc}(\Rt)$ be a weak solution of the equation (\ref{Harmonic-map}). If $\vec{\nabla} \otimes \V \in \dot{M}^{2,p}(\Rt)$, $p>3$ and $\mathbb{G}$ verifies (\ref{Cond-Th-1}) with $k\geq 0$, then we have $\V \in \mathcal{W}^{k+2,p}(\Rt)$. In addition, $\partial^\alpha \V $ is a $(1-3/p)-$H\"older continuous function for all $\vert \alpha \vert \leq k+1$. 
\end{enumerate}		
\end{Corollaire}	
\begin{Remarque} Of course, Remark \ref{Rmk} holds true for the equations (\ref{NS}) and (\ref{Harmonic-map}) in the homogeneous case  $\mathbb{F}=\mathbb{G}=0$. 	 
\end{Remarque}	

Let us briefly explain the general strategy in the proof of Theorem \ref{Th}. The proof bases on two key ideas. First, by assuming $\U, \vec{\nabla}\otimes \V \in \dot{M}^{2,p}(\Rt)$ and by using   the framework of an auxiliary parabolic system  (\ref{Sys-aux}), we prove that  $\U$ and $\vec{\nabla} \otimes \V$ are bounded functions on $\Rt$. Thereafter,  we use a bootstrap argument to show that  $\U \in \mathcal{W}^{k+2,p}(\Rt)$ and $\V \in \mathcal{W}^{k+2,p}(\Rt)$. 

\medskip 

These ideas can also  be  applied to other relevant coupled system of the fluid dynamics. This system is  the time-independent magneto-hydrodynamic equations  which describe the steady state of the magnetic properties of electrically conducting fluids, including plasma and liquid metals  \cite{Schnack}:  

\begin{equation}\label{MHD}
\left\{  \begin{array}{ll}\vspace{2mm}
-\Delta \U + \text{div}(\U \otimes \U) - \text{div}(\B \otimes \B) + \vec{\nabla} P= \text{div}(\mathbb{F}), \quad   \text{div}(\U)=0, \\ 
-\Delta \B + \text{div}(\B \otimes \U) - \text{div}(\U \otimes \B) =\text{div}(\mathbb{G}), \quad  \quad \text{div}(\B)=0,
\end{array}
\right. 
\end{equation}
Here,  $\U: \Rt \to \Rt$ and $P: \Rt \to \R$ always denote the velocity and the pressure of the fluid respectively. Moreover,  $\B: \Rt \to \Rt$ is the magnetic field. Furthermore, $\text{div}(\mathbb{F})$ and $\text{div}(\mathbb{G})$ are the external forces acting on this system.  

\medskip

Our second result states as follows:

\begin{Theoreme}\label{Th2} Let $\U \in L^{2}_{loc}(\Rt)$, $\B \in L^{2}_{loc}(\Rt)$, $P \in \mathcal{D}'(\Rt)$ be a weak solution of the system (\ref{MHD}). We assume $\U, \B \in \dot{M}^{2,p}(\Rt)$ with $p>3$. If for $k\geq 0$ the functions $\mathbb{F}$ and $\mathbb{G}$ verify (\ref{Cond-Th-1}) then we have $\U \in \mathcal{W}^{k+2,p}(\Rt)$, $\B \in \mathcal{W}^{k+2,p}(\Rt)$ and $P \in \mathcal{W}^{k+1,p}(\Rt)$. Moreover, for $\vert \alpha \vert\leq k+1$, $\partial^{\alpha} \U$ and $\partial^{\alpha}\B$ are H\"older continuous functions with  exponent $1-3/p$, while this fact holds true for $\partial^\alpha P$ with $\vert \alpha \vert \leq k$. 
\end{Theoreme}	

The system (\ref {MHD}) actually has a  simpler  structure than the system (\ref {EickLes}) as   the four nonlinear terms  have the same writing. Consequently,  they are treated similarly provided that $\U$ and $\B$ have the same properties. Thus,  we are able to adapt the ideas above  to obtain this new criterion for  weak solutions  in the setting of the Morrey spaces. 

\medskip

When $\mathbb{F}, \mathbb{G} \in L^2(\Rt)$, and consequently the external forces verify $\text{div}(\mathbb{F}), \text{div}(\mathbb{G}) \in \dot{H}^{-1}(\Rt)$, it is well-known that the system (\ref {MHD})  has \emph{finite energy} weak solutions $\U, \, \B \in \dot{H}^{1}(\Rt)$ and $P \in \dot{H}^{1/2}(\Rt)+ L^{2}(\Rt)$ such that $\Vert \U \Vert^{2}_{\dot{H}^1}+ \Vert \B \Vert^{2}_{\dot{H}^1} \leq c \Vert  \text{div}(\mathbb{F}) \Vert^{2}_{\dot{H}^{-1}} + c \Vert  \text{div}(\mathbb{G}) \Vert^{2}_{\dot{H}^{-1}}$.  See, for instance, the Theorem $16.2$ of the book \cite{PGLR1} for a proof in the case of the stationary Navier-Stokes equations (\ref{NS}), which can be easily adapted to the system (\ref {MHD}).  

\medskip 

By recalling that we have the following embedding:  $\dot{H}^{1}(\Rt) \subset L^{6}(\Rt) \subset \dot{M}^{2,6}(\Rt)$, the following corollary gives us a regularity criterion for the finite energy weak solutions of (\ref {MHD}).

\begin{Corollaire} Let $\mathbb{F}, \mathbb{G} \in L^2(\Rt)$ and let $(\U,\B)\in \dot{H}^{1}(\Rt)$ be a weak solution of the system  (\ref {MHD}). If  for $k\geq 0$ the functions $\mathbb{F}$ and $\mathbb{G}$ verify (\ref{Cond-Th-1}) then we have $\U \in \mathcal{W}^{k+2,6}(\Rt)$, $\B \in \mathcal{W}^{k+2,6}(\Rt)$ and $P \in \mathcal{W}^{k+1,6}(\Rt)$. Moreover, if $\mathbb{F}= \mathbb{G}=0$ we have $(\U, \B, P) \in \mathcal{C}^{\infty}(\Rt)$. 
\end{Corollaire}	

For these  finite energy weak solutions  we are able to go further in the study of their regularity. We recall that for a parameter $b >0$ we define the weighted exponential operator  $\ds{e^{b \, \sqrt{-\Delta}}}$ as $\ds{\mathcal{F} \Big(e^{b \, \sqrt{-\Delta}}\, \varphi \Big) (\xi)= e^{b\, \vert \xi \vert} \widehat{\varphi} (\xi)}$, for $\varphi \in \mathcal{S}(\Rt)$.   Thereafter, for a parameter $s  \in \R$ we define the Gevrey class 
\[ G^{s}_{b} (\Rt)= \left\{  f \in \dot{H}^{s}(\Rt): \,  e^{b\, \sqrt{-\Delta}} \, f \in \dot{H}^s(\Rt) \right\}.\]
For  $\vert s \vert < 3/2$, $\ds{G^{s}_{b} (\Rt)}$  is a Banach space with the norm  $\Vert e^{b \, \sqrt{-\Delta}} (\cdot) \Vert_{\dot{H}^s}$. Moreover, for $s\geq 0$  the functions in  $\ds{G^{s}_{b} (\Rt)}$ are analytic. Thus, our third result writes down as follows: 
\begin{Theoreme}\label{Th3} Let $\mathbb{F}, \mathbb{G} \in L^{2}(\Rt)$. For $b>0$ we assume $\mathbb{F},\mathbb{G} \in G^{0}_{b}(\Rt)$. Then, there exists $b_1>0$ such that all  the finite energy weak solutions $(\U, \B) \in \dot{H}^{1}(\Rt)$ and $P \in \dot{H}^{1/2}(\Rt) + L^2(\Rt)$ of the system  (\ref {MHD}) associated to $\mathbb{F}$ and $\mathbb{G}$ verify $\U \in G^{1}_{b_1}(\Rt)$, $\B \in G^{1}_{b1}(\Rt)$ and $P \in G^{1/2}_{b1}(\Rt)+ G^{0}_{b}(\Rt)$.  
\end{Theoreme}	
Consequently, we obtain that  $\U, P$ and 	$\B$ are analytic functions and they admit  holomorphic extensions to the strip $\ds{ \left\{  (x+iz) \in \mathbb{C}^3 : \, \vert z \vert < \min(b, b_1) \right\}}$. Moreover, in the homogeneous case $\mathbb{F}=\mathbb{G}=0$ we have:
\begin{Corollaire}\label{Col3} All finite energy solutions $\U, \B \in \dot{H}^{1}(\Rt)$ and $P\in \dot{H}^{1/2}(\Rt)$ of the homogeneous (MHD) system verify  $\U, \B \in G^{1}_{b_1}(\Rt)$ and $P \in G^{1/2}_{b_1}(\Rt)$,  for any $b_1>0$. 
\end{Corollaire}	

To close this section, let us mention that both results given in Theorems \ref{Th2} and \ref{Th3} respectively could be  also adapted to prove analogous theorems for other coupled systems of the fluid dynamics. For instance, the stationary Boussinesq system: 
\begin{equation*}
\begin{cases}
-\Delta \U + \text{div}(\U \otimes \U) + \vec{\nabla} P = \theta \vec{e}_3 + \text{div}(\mathbb{F}), \quad \text{div}(\U)=0, \\
-\Delta \theta + \text{div}(\theta\, \U)= \text{div}(\mathbb{G}).
\end{cases} 
\end{equation*}
where,  the unknown $\theta : \Rt \to \R$ is the temperature of the fluid and $\vec{e}_3 $ denotes the third vector of the canonical basis in $\Rt$. On the other hand, the stationary tropical climate model: 
\begin{equation*}
\begin{cases}
-\Delta \U + \text{div}(\U \otimes \U )  + \text{div}(\vec{D} \otimes \vec{D} ) + \vec{\nabla} P = \text{div}(\mathbb{F}), \\
-\Delta \vec{D} + \text{div}(\U \otimes \vec{D}) + \text{div}(\vec{D} \otimes \U) + \vec{\nabla} \theta=\text{div}(\mathbb{G}), \\
-\Delta \theta +  \U \cdot \nabla \theta + \text{div}(\vec{D})=0,  \quad \text{div}(\U)=0,
\end{cases} 
\end{equation*}
where $\vec{D}:\Rt \to \Rt$ stands for the  baroclinic mode of the velocity field $\U$. 
	 
\section{Some well-known results}\label{Sec:Well-Known results}
For the reader's convenience,  we  summarize  here some  well-known results which will be useful in the sequel.  For $1<r<p$ and $1<p<+\infty$, we consider the homogeneous Morrey space $\dot{M}^{r,p}(\Rt)$, which is defined as in (\ref{Def-Morrey}) with $r$ instead of $2$. 

\begin{Lemme}[Page  $169$ of  \cite{PGLR1}]\label{Lem-U1}  The space  $\dot{M}^{r,p}(\Rt)$   is stable under convolution with functions in the space $L^1(\Rt)$ and we have  $\ds{\Vert g \ast f \Vert_{\dot{M}^{r,p}} \leq c \Vert g \Vert_{L^1} \Vert f \Vert_{\dot{M}^{r,p}}}$. 
\end{Lemme}	
\begin{Lemme}\label{Lem-U5} Let $f \in \dot{M}^{r,p}(\Rt) \cap L^{\infty}(\Rt)$. Then, for all $1\leq \sigma < +\infty$ we have $f \in \dot{M}^{ r\sigma, p\sigma}(\Rt)$ and the following estimate holds:  $\ds{\Vert f \Vert_{\dot{M}^{r\sigma, p\sigma}}\leq c\, \Vert f \Vert^{\frac{1}{\sigma}}_{\dot{M}^{r,p}}\, \Vert f \Vert^{1-\frac{1}{\sigma}}_{L^{\infty}}}$. 
\end{Lemme}	
\pv  By the interpolation inequalities with parameter $\frac{1}{\sigma}$ we write: 
\[ \left( \int_{B(x_0,R)} \vert f(x) \vert^{r\sigma} dx\right)^{\frac{1}{r\sigma}} \leq c \, \left(\int_{B(x_0,R)} \vert f(x)\vert^{r} dx \right)^{\frac{1}{r \sigma }}\, \Vert f \Vert^{1-\frac{1}{\sigma}}_{L^{\infty}},\]
then, we multiply each side by $R^{\frac{3}{p\sigma}-\frac{3}{r\sigma}}$ to obtain:  
 \[ R^{\frac{3}{p\sigma}-\frac{3}{r\sigma}} \left( \int_{B(x_0,R)} \vert f(x) \vert^{r\sigma} dx\right)^{\frac{1}{r\sigma}}  \leq c\, \left[  R^{\frac{3}{p} - \frac{3}{r}} \, \left(\int_{B(x_0,R)} \vert f(x)\vert^{r} dx \right)^{\frac{1}{r}}\, \right]^{\frac{1}{\sigma}} \, \Vert f \Vert^{1-\frac{1}{\sigma}}_{L^{\infty}}.\]
 Finally, as $\vert B(x_0, R) \vert \simeq R^3$ we obtain $\Vert f \Vert_{\dot{M}^{r\sigma, p\sigma}}\leq c\, \Vert f \Vert^{\frac{1}{\sigma}}_{\dot{M}^{r,p}}\, \Vert f \Vert^{1-\frac{1}{\sigma}}_{L^{\infty}}$.\finpv
\begin{Lemme}[Page $171$ of  \cite{PGLR1}]\label{Lem-U2} Let $t>0$ and let $h_t$ be the heat kernel. The following estimate holds  $\ds{t^{\frac{3}{2p}} \Vert h_t \ast f \Vert_{L^{\infty}} \leq c \Vert f \Vert_{\dot{M}^{r,p}}}$.
\end{Lemme} 	
This estimate is a direct consequence of the continuous embedding $\dot{M}^{r,p}(\Rt) \subset \dot{B}^{-\frac{3}{2}, \infty}_{\infty}(\Rt)$. We recall that the homogeneous Besov space $\dot{B}^{-\frac{3}{2}, \infty}_{\infty}(\Rt)$ can be characterized  as the space of temperate distributions $f\in \mathcal{S}'(\Rt)$ such that $\ds{\sup_{t>0}\, t^{\frac{3}{2 p}}\Vert h_t \ast f \Vert_{L^{\infty}} <+\infty}$.
\begin{Lemme}[Lemme $4.2$ of \cite{Kato}]\label{Lem-U3} For $i=1,2,3$ let  $\ds{\mathcal{R}_i = \frac{\partial_i}{\sqrt{-\Delta}}}$ be the i-th Riesz transform. Then, for $i,j=1,2,3$ the operator $\mathcal{R}_{i}\mathcal{R}_{j}$ is continuous in the space $\dot{M}^{r,p}(\Rt)$ and we have $\ds{\Vert \mathcal{R}_i \mathcal{R}_j (f) \Vert_{\dot{M}^{r,p}} \leq c \Vert f  \Vert_{\dot{M}^{r,p}}}$. 
\end{Lemme}	 
Finally, we shall use the following result linking the Morrey spaces and the H\"older regularity of functions. 
\begin{Lemme}[Proposition $3.4$ of \cite{GigaMiyakawa}]\label{Lem-U4} Let $f \in \mathcal{S}'(\Rt)$ such that $\vec{\nabla} f \in \dot{M}^{1,p}(\Rt)$, with $p>3$. There exists a constant $C>0$ such that for all  $x,y\in \Rt$ we have $\ds{\vert f(x)-f(y) \vert \leq C\, \Vert \vec{\nabla} f \Vert_{\dot{M}^{1,p}} \, \vert x-y \vert^{1-3/p}}$. 
\end{Lemme} 
Recall that the Morrey  space $\dot{M}^{1,p}(\Rt)$ is defined as the space of locally finite Borel measures $d\mu$ such that 
\[ \sup_{x_0 \in \Rt,\, R>0} R^{\frac{3}{p}} \left( \frac{1}{\vert B(x_0, R)\vert} \int_{B(x_0,R)} d\vert \mu \vert (x) \right)<+\infty. \]	  

\section{Proof of the Theorem \ref{Th}}
For the sake of clearness, we shall divide the proof  in three main steps.

\medskip

{\bf Step 1. The auxiliary parabolic system}.  Our starting point is the study of the following auxiliary parabolic system. Let  $\V: \Rt \to \mathbb{S}^2$ be the vector field  given in  Definition \ref{Def-weakSol-EL}. Moreover,   $\P$ stands for the the Leray's projector.  We consider the initial value problem for the parabolic coupled system: 
\begin{equation}\label{Sys-aux}
\left\{ \begin{array}{ll}\vspace{2mm} 
\partial_t \vu 	-\Delta \vu + \P(\text{div} (\vu\otimes \vu))   + \P(\text{div} ( \bV  \odot \bV))= \P\left( \text{div}(\mathbb{F})\right), \quad \text{div}(\vu)=0, \\ \vspace{2mm} 
\partial_t \bV 	-\Delta  \bV  + \vec{\nabla} \otimes  (\vu \,  \bV )  - \vec{\nabla} \otimes ( \vert \bV \vert^2\, \V) = \vec{\nabla} \otimes \left(\text{div}(\mathbb{G})\right),  \\ 
\vu(0,\cdot)=\vu_0,\, \quad \bV (0,\cdot)= \bV_0,
\end{array} \right. 
\end{equation}
where, the   vector field $\vu=(u_1, u_2, u_3)$  and the matrix $\bV = (v_{i,j})_{1\leq i,j \leq 3}$   are the unknowns.   We emphasize that in the second equation the vector field $\V$ is given.

\medskip 

 For a time $0<T<+\infty$, we denote  $\mathcal{C}_{*}([0,T], \dot{M}^{2,p}(\Rt))$ the functional space of bounded and  weak$-*$ continuous  functions from $[0,T]$ with values in the Morrey space $\dot{M}^{2,p}(\Rt)$. We prove now the following:  

\begin{Proposition}\label{Prop1} Consider the initial value problem (\ref{Sys-aux}) where $\mathbb{F}$ and $\mathbb{G}$ verify (\ref{Cond-Th-1}). Let $p>3$ and  let   $\vu_0 \in  \dot{M}^{2,p}(\Rt)$ and $ \bV_0 \in \dot{M}^{2,p}(\Rt)$  be the initial data. There exists a time $T_0>0$,  depending on $\vu_0$, $\bV_0$, $\mathbb{F}$ and $\mathbb{G}$; and there exist  $(\vu, \bV) \in \mathcal{C}_{*}([0,T_0], \dot{M}^{2,p}(\Rt))$, which is the unique solution of (\ref{Sys-aux}).  Moreover this solution verifies:   
	\begin{equation}\label{Cond-Sup}
	\sup_{0<t<T_0} t^{\frac{3}{2 p}} \Big(\Vert  \vu(t, \cdot) \Vert_{L^{\infty}}+ \Vert  \bV(t, \cdot) \Vert_{L^{\infty}}\Big) <+\infty. 
	\end{equation}
\end{Proposition} 	
\pv   Mild  solutions of the system (\ref{Sys-aux}) write down as the integral formulation:
\begin{equation}\label{equ}
\begin{split}
\vu(t,\cdot)= & \,\, e^{t \Delta}\vu_0  + \int_{0}^{t} e^{(t-s)\Delta} \,  \P \left(\text{div}(\mathbb{F}) \right)ds  + \underbrace{ \int_{0}^{t} e^{(t-s)\Delta }\P (div (\vu \otimes \vu ))(s,\cdot) ds}_{B_1(\vu,\vu)} \\
& + \underbrace{\int_{0}^{t} e^{(t-s)\Delta} \P(\text{div} (  \bV\odot \bV))(s,\cdot) ds}_{B_2(\bV, \bV)},
\end{split}
\end{equation}
and 
\begin{equation}\label{eqv}
\begin{split}
\bV(t,\cdot)= & \,\, e^{t\Delta} \bV_0 + \int_{0}^{t}e^{(t-s)\Delta} \, \vec{\nabla} \otimes \left(\text{div}(\mathbb{G}) \right)ds  + \underbrace{ \int_{0}^{t} e^{(t-s)\Delta} \vec{\nabla} \otimes  (\vu\, \bV )(s,\cdot) ds}_{B_3(\vu, \bV)}\\
&- \underbrace{ \int_{0}^{t} e^{(t-s)\Delta}  \vec{\nabla} \otimes ( \vert  \bV \vert^2\, \V)(s,\cdot)ds}_{B_4(\bV, \bV)}.
\end{split}
\end{equation}
By the  Picard's fixed point argument,  we will solve both problems  (\ref{equ}) and (\ref{eqv}) in the Banach space 
\[E_T = \left\{  f \in \mathcal{C}_{*}([0,T], \dot{M}^{2,p}(\Rt)): \sup_{0<t<T} t^{\frac{3}{2 p}} \Vert f(t,\cdot)\Vert_{L^{\infty}}<+\infty \right\}, \]
 with the norm 
\[ \Vert f \Vert_{E_T}= \sup_{0\leq t \leq T} \Vert f(t,\cdot)\Vert_{\dot{M}^{2,p}}+\sup_{0<t<T} t^{\frac{3}{2p}} \Vert f(t,\cdot)\Vert_{L^{\infty}}. \] 
Let us mention that for $f_1, f_2\in E_T$, for the sake of simplicity, we shall write $\Vert (f_1, f_2) \Vert_{E_T}= \Vert f_1\Vert_{E_T}+ \Vert f_2 \Vert_{E_T}$. 

\medskip

We start by studying the  linear terms in (\ref{equ}) and (\ref{eqv}). As $\vu_0 \in \dot{M}^{2,p}(\Rt)$ and $\bV_0 \in \dot{M}^{2,p}(\Rt)$ by Lemma \ref{Lem-U1}  we have $\ds{\Vert (e^{t\Delta} \vu_0, e^{t\Delta} \bV_0) \Vert_{\dot{M}^{2,p}} \leq c \Vert (\vu_0, \bV_0) \Vert_{\dot{M}^{2,p}}}$, hence we obtain  $\ds{e^{t \Delta}\vu_0 \in \mathcal{C}_{*}([0,T], \dot{M}^{2,p}(\Rt))}$ and $\ds{e^{t \Delta}\bV_0 \in \mathcal{C}_{*}([0,T], \dot{M}^{2,p}(\Rt))}$. On the other hand,  by Lemma \ref{Lem-U2} the following estimate directly follows $\ds{\sup_{0<t<T}t^{\frac{3}{2p}} \left\Vert \left(e^{t\Delta} \vu_0, e^{t\Delta} \bV_0\right)\right\Vert_{L^{\infty}}  \leq c \Vert (\vu_0, \bV_0) \Vert_{\dot{M}^{2,p}}}$. Thus, we have $e^{t \Delta}\vu_0 \in E_T$ and $e^{t\Delta} \bV_0 \in E_T$, and moreover, the following estimate holds:
\begin{equation}\label{Lin}
\left\Vert \left(e^{t\Delta} \vu_0, e^{t\Delta} \bV_0\right) \right\Vert_{E_T} \leq c \Vert (\vu_0, \bV_0) \Vert_{\dot{M}^{2,p}}. 
\end{equation} 
Thereafter, as $\mathbb{F}, \mathbb{G}$ are time independent tensors, and moreover, as we assume (\ref{Cond-Th-1}), we write:
\begin{equation*}
\begin{split}
\left\Vert \int_{0}^{t} e^{(t-s)\Delta} \left( \text{div}(\mathbb{F}), \text{div}(\mathbb{G})\right) ds  \right\Vert_{\dot{M}^{2,p}}  \leq  & \,\, \int_{0}^{t} \left\Vert e^{(t-s)\Delta} \left(\text{div}(\mathbb{F}), \text{div}(\mathbb{G}) \right) \right\Vert_{\dot{M}^{2,p}}\, ds \\
\leq & \,\,  c \left\Vert \left( \text{div}(\mathbb{F}), \text{div}(\mathbb{G})\right) \right\Vert_{\dot{M}^{2,p}} \left( \int_{0}^{t} ds\right),
\end{split}
\end{equation*} to get
\begin{equation}\label{estim-forces-1}
\sup_{0 \leq t \leq T} \left\Vert \int_{0}^{t} e^{(t-s)\Delta} \left(\text{div}(\mathbb{F}), \text{div}(\mathbb{G}) \right) ds  \right\Vert_{\dot{M}^{2,p}}  \leq c \, T\, \left\Vert \left(\text{div}(\mathbb{F}),\text{div}(\mathbb{G}) \right) \right\Vert_{\dot{M}^{2,p}}.
\end{equation}
On the other hand,  we remark that by Lemma \ref{Lem-U2} we have
\[ \left\Vert e^{(t-s)\Delta}\left(\text{div}(\mathbb{F}),\text{div}(\mathbb{G})\right) \right\Vert_{L^{\infty}} \leq c\, (t-s)^{-\frac{3}{2p}} \left\Vert \left(\text{div}(\mathbb{F}),\text{div}(\mathbb{G})\right) \right\Vert_{\dot{M}^{2,p}},\]
and then we can  write
\begin{equation*}
\begin{split}
t^{\frac{3}{2 p}}\, \left\Vert \int_{0}^{t} e^{(t-s)\Delta}\left(\text{div}(\mathbb{F}), \text{div}(\mathbb{G})\right) ds \right\Vert_{L^\infty} & \leq t^{\frac{3}{2 p}} \, \int_{0}^{t} \left\Vert e^{(t-s)\Delta}\left(\text{div}(\mathbb{F}),\text{div}(\mathbb{G})\right) \right\Vert_{L^{\infty}} ds \\
& \leq c\, t^{\frac{3}{2 p}} \, \int_{0}^{t} (t-s)^{-\frac{3}{2p}} \left\Vert \left(\text{div}(\mathbb{F}),\text{div}(\mathbb{G})\right) \right\Vert_{\dot{M}^{2,p}} ds\\
& \leq c\,  t^{\frac{3}{2 p}}  \,  \left\Vert \left( \text{div}(\mathbb{F}), \text{div}(\mathbb{G})\right) \right\Vert_{\dot{M}^{2,p}}  \, \left( \int_{0}^{t} (t-s)^{-\frac{3}{2p}} \, ds  \right) \\
&\leq c \, t\, \left\Vert \left( \text{div}(\mathbb{F}), \text{div}(\mathbb{G})\right) \right\Vert_{\dot{M}^{2,p}}. 
\end{split}
\end{equation*}
We thus obtain
\begin{equation}\label{estim-forces-2}
\sup_{0 <t < T} t^{\frac{3}{2 p}} \, \left\Vert \int_{0}^{t} e^{(t-s)\Delta}\left(\text{div}(\mathbb{F}),\text{div}(\mathbb{G})\right) ds \right\Vert_{L^\infty}  \leq c \, T\, \left\Vert \left(\text{div}(\mathbb{F}),\text{div}(\mathbb{G}) \right) \right\Vert_{\dot{M}^{2,p}}.
\end{equation}
By the estimates  (\ref{estim-forces-1}) ans (\ref{estim-forces-2}) we get
\begin{equation}\label{estim-forces}
\left\Vert \int_{0}^{t} e^{(t-s)\Delta}\left(\text{div}(\mathbb{F}), \text{div}(\mathbb{G}) \right) ds \right\Vert_{E_T}  \leq c\, T\, \left\Vert \left(\text{div}(\mathbb{F}), \text{div}(\mathbb{G}) \right) \right\Vert_{\dot{M}^{2,p}}.
\end{equation}

We study now the bilinear terms in (\ref{equ}) and (\ref{eqv}). First,  the terms  $B_{1}(\vu,\vu)$ and $B_2(\bV, \bV)$ in (\ref{equ}) are estimated as follows: 
\begin{equation}\label{Estim-non-lin-1}
\sup_{0\leq t \leq T} \left\Vert B_1(\vu, \vu) + B_2(\bV, \bV) \right\Vert_{\dot{M}^{2,p}}  \leq  c\, T^{\frac{1}{2}-\frac{3}{2p}}\,  \left\Vert \left(\vu, \bV\right) \right\Vert^{2}_{E_T},
\end{equation}
where, as $p>3$ then we have  $\frac{1}{2}-\frac{3}{2p}>0$. Indeed, by Lemma \ref{Lem-U1}, by the  well-known estimate on the heat kernel: $\ds{\Vert \vec{\nabla}h_{(t-s)} (\cdot)\Vert_{L^1} \leq \frac{c}{(t-s)^{1/2}}}$, and moreover, by Lemma \ref{Lem-U3} (hence $\P$ is continuous in $\dot{M}^{2,p}(\Rt)$) we have: 
\begin{equation}\label{BiLin1}
\begin{split}
& \sup_{0\leq t \leq T} \left\Vert B_1(\vu, \vu) + B_2(\bV, \bV) \right\Vert_{\dot{M}^{2,p}} \\
=  &  \,\, \sup_{0\leq t \leq T} \left\Vert  \int_{0}^{t} e^{(t-s)\Delta }\P (div (\vu \otimes \vu ))(s,\cdot) ds + \int_{0}^{t} e^{(t-s)\Delta} \P(\text{div} ( \bV \odot \bV))(s,\cdot) ds \right\Vert_{\dot{M}^{2,p}}\\
\leq &\,\,  c\,  \sup_{0\leq t \leq T} \int_{0}^{t}  \left\Vert   e^{(t-s)\Delta } (div (\vu \otimes \vu ))(s,\cdot)+  e^{(t-s)\Delta} (\text{div} ( \bV \odot \bV))(s,\cdot)\right\Vert_{\dot{M}^{2,p}}  ds\\
\leq & \,\,c\,  \sup_{0\leq t \leq T} \int_{0}^{t} \frac{1}{(t-s)^{1/2}} \left( \Vert \vu(s,\cdot) \otimes  \vu(s,\cdot)\Vert_{\dot{M}^{2,p}}+\Vert \bV(s,\cdot) \odot \bV(s,\cdot) \Vert_{\dot{M}^{2,p}}\right) ds\\
\leq &\,\, c\, \sup_{0\leq t \leq T} \int_{0}^{t} \frac{1}{(t-s)^{\frac{1}{2}}\,s^{\frac{3}{2 p}} }  \left( (s^{\frac{3}{2p}}\Vert \vu(s,\cdot) \Vert_{L^{\infty}}) \Vert \vu(s,\cdot)\Vert_{\dot{M}^{2,p}} +  (s^{\frac{3}{2p}}\Vert \bV(s,\cdot) \Vert_{L^{\infty}}) \Vert \bV(s,\cdot)\Vert_{\dot{M}^{2,p}} \right)ds\\
\leq & \,\,c\, T^{\frac{1}{2}-\frac{3}{2p}}\,  \left\Vert \left(\vu, \bV\right) \right\Vert^{2}_{E_T}. 
\end{split}
\end{equation}  

We will prove now  the following estimate:
\begin{equation}\label{Estim-non-lin-2}
\sup_{0\leq t \leq T} t^{\frac{3}{2p}} \left\Vert B_1(\vu, \vu) + B_2(\bV, \bV) \right\Vert_{L^{\infty}} \leq  c\, T^{\frac{1}{2}-\frac{3}{2p}}\,  \left\Vert \left(\vu, \bV\right) \right\Vert^{2}_{E_T}.
\end{equation}
We write: 
\begin{equation*}
\begin{split}
&\,\,  \sup_{0\leq t \leq T} t^{\frac{3}{2p}} \left\Vert B_1(\vu, \vu) + B_2(\bV, \bV) \right\Vert_{L^{\infty}} \\
= & \,\, \sup_{0<t<T} t^{\frac{3}{2p}}  	 \left\Vert  \int_{0}^{t} e^{(t-s)\Delta }\P (div (\vu \otimes \vu ))(s,\cdot) ds + \int_{0}^{t} e^{(t-s)\Delta} \P(\text{div} ( \bV\odot \bV))(s,\cdot) ds \right\Vert_{L^{\infty}} \\
\leq & \,\,  \sup_{0<t<T} t^{\frac{3}{2p}}    \int_{0}^{t}  \left\Vert  e^{(t-s)\Delta } \P \, \text{div} \left( \vu \otimes \vu  + \bV \odot \bV \right)(s,\cdot)  \right\Vert_{L^{\infty}} ds=(a).
\end{split}
\end{equation*}
Here, we recall that the operator $e^{(t-s)\Delta} \P(div(\cdot))$ writes down as  a matrix of convolution operators (in the spatial variable) whose kernels $K_{i,j}$ verify $\ds{\vert K_{i,j}(t-s, x )\vert \leq \frac{c}{((t-s)^{1/2}+\vert x \vert)^{4}}}$, see Proposition $11.1$ of \cite{PGLR}.  Then, we have  $\ds{\Vert K_{i,j}(t-s,\cdot) \Vert_{L^1} \leq \frac{c}{(t-s)^{1/2}}}$; and we can write:
\begin{equation*}
\begin{split}
(a)\leq &  \,\, c\,  \sup_{0\leq t \leq T} t^{\frac{3}{2p}}   \int_{0}^{t} \frac{1}{(t-s)^{1/2}} \left( \Vert \vu(s,\cdot) \otimes  \vu(s,\cdot)\Vert_{L^{\infty}}+\Vert \bV(s,\cdot) \odot \bV(s,\cdot) \Vert_{L^{\infty}}\right) ds\\
\leq & \,\,  c\,  \sup_{0\leq t \leq T} t^{\frac{3}{2p}}   \int_{0}^{t} \frac{ds}{(t-s)^{1/2} s^{\frac{3}{p}}} \left(\left( s^{\frac{3}{2p}} \Vert \vu(s,\cdot) \Vert_{L^{\infty}}\right)^2+\left( s^{\frac{3}{2p}}\Vert\bV(s,\cdot)  \Vert_{L^{\infty}}\right)^2\right) ds\\
\leq & \,\,  c\, \left( \sup_{0\leq t \leq T} t^{\frac{3}{2p}}    \int_{0}^{t} \frac{ds}{(t-s)^{1/2} s^{\frac{n}{p}}} \right) \left\Vert \left(\vu, \bV\right) \right\Vert^{2}_{E_T}. \\
\leq & \,\,   \left(c\,  \sup_{0\leq t \leq T}\left[ t^{\frac{3}{2p}}   \int_{0}^{t/2} \frac{ds}{(t-s)^{1/2} s^{\frac{3}{p}}} +t^{\frac{3}{2p}}   \int_{t/2}^{t} \frac{ds}{(t-s)^{1/2} s^{\frac{3}{p}}} \right]\right)  \left\Vert \left(\vu, \bV\right) \right\Vert^{2}_{E_T}\\
\leq &  \,\,  c\,\left(  \sup_{0\leq t \leq T}\left[ t^{\frac{3}{2p}-\frac{1}{2}}    \int_{0}^{t/2} \frac{ds}{s^{3/p}} +t^{\frac{3}{2p}-\frac{3}{p}}  \int_{t/2}^{t} \frac{ds}{(t-s)^{1/2} } \right]\right)  \left\Vert \left(\vu, \bV\right) \right\Vert^{2}_{E_T}\\
\leq &  \,\,  c\, T^{\frac{1}{2}-\frac{3}{2p}}\,  \left\Vert \left(\vu, \bV\right) \right\Vert^{2}_{E_T}.
\end{split}
\end{equation*}

We study now the bilinear terms  $B_3(\vu, \bV)$ and $B_4(\bV, \bV)$ in (\ref{eqv}). Precisely, we shall prove the estimates: 
\begin{equation}\label{Estim-non-lin-3}
\sup_{0\leq t \leq T} \left\Vert B_3(\vu, \bV) +  B_4(\bV, \bV) \right\Vert_{\dot{M}^{2,p}}  \leq  c\, T^{\frac{1}{2}-\frac{3}{2p}} \, \left\Vert \left( \vu, \bV \right) \right\Vert^{2}_{E_T},
\end{equation}
and
\begin{equation}\label{Estim-non-lin-4}
\sup_{0\leq t \leq T}  t^{\frac{3}{2 p}} \left\Vert B_{3}(\vu, \bV)+ B_4(\bV,\bV) \right\Vert_{L^{\infty}}  \leq  c\, T^{\frac{1}{2}-\frac{3}{2p}} \, \left\Vert \left( \vu, \bV \right) \right\Vert^{2}_{E_T}. 
\end{equation}
In order to prove the estimate (\ref{Estim-non-lin-3}), we recall that the vector field $\V$ verifies $\V(x)=1$ for \emph{a.e.} $x \in \Rt$ (see Definition \ref{Def-weakSol-EL}) and then we have $\Vert \V\Vert_{L^{\infty}}=1$. With this information at hand,  for the first in the norm $\Vert \cdot \Vert_{E_T}$ we are able to write the following estimates: 
\begin{equation}\label{BiLin3}
\begin{split}
& \sup_{0\leq t \leq T} \left\Vert B_3(\vu, \bV) +  B_4(\bV, \bV) \right\Vert_{\dot{M}^{2,p}} \\
=&  \,\, \sup_{0\leq t \leq T} \left\Vert 	 \int_{0}^{t} e^{(t-s)\Delta} \vec{\nabla} \otimes  (\vu \, \bV)(s,\cdot) ds- \int_{0}^{t} e^{(t-s)\Delta}  \vec{\nabla} \otimes ( \vert  \bV \vert^2\, \V)(s,\cdot)ds \right\Vert_{\dot{M}^{2,p}}\\
\leq & \,\, c\, \sup_{0\leq t \leq T} \, \int_{0}^{t} \frac{1}{(t-s)^{1/2}} \left( \Vert  \vu \, \bV(s,\cdot) \Vert_{\dot{M}^{2,p}}+ \Vert \vert \bV(s,\cdot) \vert^2 \V \Vert_{\dot{M}^{2,p}}  \right) ds \\
\leq & \,\, c\, \sup_{0\leq t \leq T}\, \int_{0}^{t} \frac{1}{(t-s)^{1/2}} \left( \Vert  \vu(s,\cdot) \Vert_{\dot{M}^{2,p}} \Vert \bV(s,\cdot) \Vert_{L^{\infty}} + \Vert \vert \bV(s,\cdot) \vert^2 \Vert_{\dot{M}^{2,p}} \Vert \V\Vert_{L^{\infty}} \right) ds\\
\leq &  \,\, c\, \sup_{0\leq t \leq T} \, \int_{0}^{t} \frac{1}{(t-s)^{1/2}} \left( \Vert \vu (s,\cdot) \Vert_{\dot{M}^{2,p}} \Vert \bV(s,\cdot)\Vert_{L^{\infty}}  + \Vert \vert \bV(s,\cdot) \vert^2 \Vert_{\dot{M}^{2,p}} \right) ds\\
\leq &  \,\, c\, \sup_{0\leq t \leq T}\, \int_{0}^{t} \frac{1}{(t-s)^{1/2}} \left( \Vert \vu (s,\cdot) \Vert_{\dot{M}^{2,p}} \Vert \bV(s,\cdot)\Vert_{L^{\infty}}  + \Vert  \bV(s,\cdot)  \Vert_{\dot{M}^{2,p}} \Vert \bV(s,\cdot) \Vert_{L^{\infty}} \right) ds\\
\leq &  \,\, c\, \left[ \sup_{0\leq t \leq T}\, \int_{0}^{t} \frac{ds}{(t-s)^{1/2} s^{3/2p}}\right] \left\Vert \left( \vu, \bV \right) \right\Vert^{2}_{E_T} \leq  c\, T^{\frac{1}{2}-\frac{3}{2p}} \, \left\Vert \left( \vu, \bV \right) \right\Vert^{2}_{E_T}.
\end{split}
\end{equation}
In order to prove the estimate (\ref{Estim-non-lin-4}), we essentially follow the estimates performed in (\ref{Estim-non-lin-2}) to obtain: 
\begin{equation}\label{BiLin4}
\begin{split}
& \sup_{0\leq t \leq T}  t^{\frac{3}{2 p}} \left\Vert B_{3}(\vu, \bV)+ B_4(\bV,\bV) \right\Vert_{L^{\infty}} \\
=& \,\, \sup_{0\leq t \leq T} t^{\frac{3}{2 p}} \left\Vert 	 \int_{0}^{t} e^{(t-s)\Delta} \vec{\nabla} \otimes  (\vu \, \bV)(s,\cdot) ds- \int_{0}^{t} e^{(t-s)\Delta}  \vec{\nabla} \otimes ( \vert  \bV \vert^2\, \V)(s,\cdot)ds \right\Vert_{L^{\infty}}\\
\leq & \,\, c\, T^{\frac{1}{2}-\frac{3}{2p}} \, \left\Vert \left( \vu, \bV \right) \right\Vert^{2}_{E_T}. 
\end{split}
\end{equation}

Summarizing, we obtain the following estimate for the four bilinear terms in the equations (\ref{equ}) and (\ref{eqv}): 
\begin{equation}\label{BiLin}
\left\Vert B_{1}(\vu,\vu) \right\Vert_{E_T}+ \left\Vert B_2(\bV, \bV)\right\Vert_{E_T}+\left\Vert B_3(\vu,\bV)\right\Vert_{E_T}+\left\Vert B_4(\bV,\bV) \right\Vert_{E_T}\leq c\, T^{\frac{1}{2}-\frac{3}{2p}} \, \left\Vert \left( \vu, \bV \right) \right\Vert^{2}_{E_T}, \quad \frac{1}{2}-\frac{3}{2p}>0.
\end{equation}
Once we have the estimates (\ref{Lin}), (\ref{estim-forces}) and (\ref{BiLin}), for  a time $0<T_0=T_0(\vu_0, \bV_0, \text{div}(\mathbb{F}), \text{div}(\mathbb{G}))<+\infty$ small enough, the   existence and uniqueness  of a solution $(\vu,\bV)$ for the equations (\ref{equ}) and (\ref{eqv}) follow from standard arguments.  Proposition  \ref{Prop1} is proven.  \finpv 

\medskip

{\bf Step 2. The global boundness of $\U$ and $\vec{\nabla} \otimes \V$}. With the help of Proposition \ref{Prop1},  we are able to prove  the following:  
\begin{Proposition}\label{Prop2}    Let  $(\U,P,\V)$ be a  weak solution of the system  (\ref{EickLes}) given in  Definition \ref{Def-weakSol-EL}. If $\U \in \dot{M}^{2,p}(\Rt)$ and $\vec{\nabla} \otimes \V \in \dot{M}^{2,p}(\Rt)$, with $p>3$, then we have $\U \in L^{\infty}(\Rt)$ and $\vec{\nabla} \otimes \V \in L^{\infty}(\Rt)$. 
\end{Proposition} 	
\pv   In the initial value  problem (\ref{Sys-aux}), we set  the initial data $\ds{(\vu_0, \bV_0)=(\U, \vec{\nabla}\otimes \V)}$. Then, by  Proposition \ref{Prop1} there exists a time $0<T_0$ and there exists  a unique  solution $\ds{(\vu, \bV)} \in \mathcal{C}_{*}([0,T], \dot{M}^{2,p}(\Rt))$ of   (\ref{Sys-aux}) arising from  $\ds{(\U, \vec{\nabla} \otimes \V)}$. 

\medskip 

On the other hand, we have the following key remark. First,  we apply the Leray's projector $\P$  in the first equation of the system (\ref{EickLes}). Newt, we apply the operator $\vec{\nabla} \otimes (\cdot)$ in the second equation of this system. Moreover, as $\U$ and $ \V$ are time-independent functions we have $\partial_t \U =0$ and $\partial_t (\vec{\nabla} \otimes \V)=0$. Thus, the  couple $\ds{(\U,\vec{\nabla} \otimes \V )}$ is also a solution of the initial value problem (\ref{Sys-aux}) with the  initial data $\ds{(\vu_0, \bV_0)=(\U, \vec{\nabla}\otimes \V)}$, and moreover, we have $\ds{(\U, \vec{\nabla}\otimes \V) \in \mathcal{C}_{*}([0,T], \dot{M}^{2,p}(\Rt))}$.

\medskip

Consequently, in the space $\ds{\mathcal{C}_{*}([0,T], \dot{M}^{2,p}(\Rt))}$ we  have   two solutions of  (\ref{Sys-aux}) with the same initial data: on the one hand, the solution $\ds{(\vu, \bV)}$ given by  Proposition \ref{Prop1} and, on the other hand,  the solution $\ds{(\U, \vec{\nabla}\otimes \V)}$. By  uniqueness we have the identity $\ds{(\vu,  \bV)= (\U, \vec{\nabla}\otimes \V)}$ and by   (\ref{Cond-Sup})  we can write 
\[  \sup_{0<t<T} t^{\frac{3}{2 p}} \left(\Vert  \U \Vert_{L^{\infty}}+ \Vert  \vec{\nabla}\otimes \V \Vert_{L^{\infty}}\right) <+\infty. \]
But, as the solution $ (\U, \vec{\nabla}\otimes \V)$ does not depend on the time variable we  have $\U \in L^{\infty}(\Rt)$ and $\vec{\nabla}\otimes \V \in L^{\infty}(\Rt)$.  Proposition  \ref{Prop2} is now proven. \finpv  

\medskip

{\bf Step 3. Estimates on high order derivatives in the Morrey spaces}. The global boundness of $\U$ and $\vec{\nabla} \otimes \V$ obtained in the previous step is the key tool to prove the following: 
\begin{Proposition}\label{Prop3}  We assume that $\mathbb{F}$ and $\mathbb{G}$ verify (\ref{Cond-Th-1}) for $k \geq 0$, and moreover, we assume $\U, \vec{\nabla} \otimes \V  \in \dot{M}^{2,p}(\Rt)$ with $p>3$. Then  we have $\U \in \mathcal{W}^{k+2,p}(\Rt)$, $\V \in \mathcal{W}^{k+2,p}(\Rt)$ and  $P \in \mathcal{W}^{k+1,p}(\Rt)$.  
\end{Proposition}	

\pv  We  will study first  the functions $\U$ and $\V$. For this,  we  get rid  (temporally) of the pressure term by applying the Leray's projector $\P$ in the first equation of the system  (\ref{EickLes}). Then, we shall consider the following coupled system:  
\begin{equation}\label{EickLesP}
\left\{ \begin{array}{ll}\vspace{2mm} 
-\Delta \U + \P ( \text{div} (\U \otimes \U) )  + \P ( \text{div} ( \vec{\nabla} \otimes  \V \odot \vec{\nabla} \otimes \V)) = \P\left(  \text{div}(\mathbb{F})\right), \\ \vspace{2mm} 
-\Delta \V + \text{div}(\V \otimes \U)   - \vert  \vec{\nabla} \otimes \V  \vert^2\, \V =  \text{div}(\mathbb{G}),  \\
\text{div}(\U)=0, \end{array} \right. 
\end{equation}
As $\U$ and $\V$ solve this system  they  verify the  following  (equivalent) integral formulations: 
\begin{equation}\label{EickLes-Int-1}
\begin{array}{ll}\vspace{2mm} 
\ds{\U = - \frac{1}{-\Delta} \left( \P ( \text{div} (\U \otimes \U) ) \right) - \frac{1}{-\Delta} \left( \P ( \text{div} ( \vec{\nabla} \otimes  \V \odot \vec{\nabla} \otimes \V)) \right)} + \frac{1}{-\Delta}\Big( \P \left(   \text{div}(\mathbb{F})\right)\Big),  
\end{array} 
\end{equation}
\begin{equation}\label{EickLes-Int-2}
\begin{array}{ll}\vspace{2mm} 
\ds{\V =- \frac{1}{-\Delta} \left( \text{div}(\V \otimes \U) \right) + \frac{1}{-\Delta} \left(  \vert  \vec{\nabla} \otimes \V  \vert^2\, \V \right)}+ \frac{1}{-\Delta} \Big(\text{div}(\mathbb{G})\Big).
 \end{array} 
\end{equation}

By using these integral formulations, we will show  that  $\partial^{\alpha} \U \in \dot{M}^{2,p}(\Rt)$  and $\partial^{\alpha} \V \in \dot{M}^{2,p}(\Rt)$ for all multi-indice $\vert \alpha \vert \leq k+2$.  We shall prove this fact by iteration respect to the order of the multi-indices $\alpha$, which  we will denote as $\vert \alpha \vert$. For the reader's convenience, in the following couple of  technical lemmas  we prove each step in the iterative argument separately.
\begin{Lemme}[The case initial case]  Recall that by Proposition \ref{Prop2} we have   $\U, \vec{\nabla} \otimes \V \in L^{\infty}(\Rt)$. Then,  for $\vert \alpha \vert \leq 2$ and for all $1\leq \sigma < +\infty$ we have $\partial^{\alpha}\U \in \dot{M}^{2\sigma, p\sigma}(\Rt)$ and $\partial^{\alpha}\V \in \dot{M}^{2\sigma,p\sigma}(\Rt)$.
\end{Lemme}	 
\pv  Due to the  coupled structure of the  equations (\ref{EickLes-Int-1}) and (\ref{EickLes-Int-2}), we  must study first the function $\V$ and then we study the function $\U$. 
\begin{enumerate}
\item[$\bullet$]  Let 	$\vert \alpha \vert=1$.  As $\vec{\nabla}\otimes \V \in \dot{M}^{2,p}(\Rt) \cap L^{\infty}(\Rt)$, by Lemma \ref{Lem-U5}   we  have $\partial^{\alpha} \V \in \dot{M}^{2\sigma, p\sigma}(\Rt)$ for all $1 \leq \sigma < +\infty$.  On the other hand, for the function $\U$, by  (\ref{EickLes-Int-1})  we have the identity: 
\begin{equation}\label{eq-der-alpha-U}
\partial^{\alpha }\U = - \frac{1}{-\Delta} \left( \P ( \partial^{\alpha} \text{div} (\U \otimes \U) ) \right) - \frac{1}{-\Delta} \left( \P (\partial^{\alpha} \text{div} ( \vec{\nabla} \otimes  \V \odot \vec{\nabla} \otimes \V)) \right)+  \frac{1}{-\Delta}\Big( \P\left(  \text{div}( \partial^{\alpha}\mathbb{F})\right)\Big),
\end{equation}
where  we shall verify that each term on the right-hand side belong to the space $\dot{M}^{2\sigma, p\sigma}(\Rt)$, for all $1\leq \sigma <+\infty$. For the first term, as $\U \in \dot{M}^{2,p}(\Rt) \cap L^{\infty}(\Rt)$ by Lemma \ref{Lem-U5} we have $\U \in \dot{M}^{r\sigma, p\sigma}(\Rt)$, for all $1\leq \sigma <+\infty$. Moreover, by the H\"older inequalities we also have $\U \otimes \U \in \dot{M}^{r\sigma, p\sigma}(\Rt)$ (for all $1\leq \sigma <+\infty$). Then,  as $\vert \alpha \vert=1$  the operator $\ds{\frac{1}{-\Delta} \left( \P ( \partial^{\alpha} \text{div}(\cdot)\right)}$ writes down as a linear combination of the Riesz transforms $\mathcal{R}_{i}\mathcal{R}_{j}$ with $i,j=1,2,3$; and  by  Lemma \ref{Lem-U3}   we obtain 
$\ds{\left\Vert  \frac{1}{-\Delta} \left( \P ( \partial^{\alpha} \text{div} (\U \otimes \U) ) \right) \right\Vert_{\dot{M}^{2\sigma, p \sigma}} \leq c\, \Vert \U \otimes \U  \Vert_{\dot{M}^{2\sigma, p \sigma}} <+\infty}$. 

\medskip

 The second term is similarly estimated by using now the information  $\vec{\nabla}\otimes \V \in \dot{M}^{r\sigma, p\sigma}(\Rt)$.  

\medskip

We study  the third term. By Lemma \ref{Lem-U3} we write  $\ds{\left\Vert \frac{1}{-\Delta}\Big( \P\left(  \text{div}( \partial^{\alpha}\mathbb{F})\right)\Big) \right\Vert_{\dot{M}^{2\sigma, p \sigma}} \leq c\, \Vert \mathbb{F} \Vert_{\dot{M}^{2\sigma, p\sigma}}}$.
Then,  by (\ref{Cond-Th-1}) we have $\mathbb{F}\in \dot{M}^{2,p}\cap L^{\infty}(\Rt)$, and by Lemma \ref{Lem-U5}  we get $\ds{\Vert \mathbb{F} \Vert_{\dot{M}^{2\sigma, p\sigma}} \leq c\, \Vert \mathbb{F} \Vert^{\frac{1}{\sigma}}_{\dot{M}^{2,p}}\, \Vert \mathbb{F} \Vert^{1-\frac{1}{\sigma}}_{L^{\infty}} <+\infty}$. 

\item[$\bullet$] Let $\vert \alpha \vert =2$.   By (\ref{EickLes-Int-2}) we write 
\begin{equation}\label{eq-der-V}
\partial^{\alpha}\V =- \frac{1}{-\Delta} \left(  \partial^{\alpha} \text{div}(\V \otimes \U) \right) + \frac{1}{-\Delta} \left( \partial^{\alpha} \left( \vert  \vec{\nabla} \otimes \V  \vert^2\, \V  \right)\right)+ \frac{1}{-\Delta} \Big(\partial^{\alpha} \text{div}(\mathbb{G})\Big).
\end{equation}
As before, we will prove that each term on the right-hand side belong to the space $\dot{M}^{2\sigma, p\sigma}(\Rt)$.  For the first term,  for $i,j=1,2,3$ we write $\partial_i (V_i U_j)= \partial_i V_i\, U_j + V_i \partial_i U_j$.  By recalling that  $\vec{\nabla} \otimes \V \in \dot{M}^{2\sigma, p\sigma}(\Rt)$,  $\U \in L^{\infty}(\Rt)$,  and moreover, by recalling that $\V \in L^{\infty}(\Rt)$ and $\partial_i \U \in \dot{M}^{2\sigma, p\sigma}(\Rt)$, 
we directly have $\ds{\text{div}(\V \otimes \U) \in \dot{M}^{2\sigma, p\sigma} (\Rt)}$. Thereafter, as $\vert \alpha\vert=2$, by Lemma \ref{Lem-U3}  the operator $- \frac{1}{-\Delta} \left(  \partial^{\alpha} (\cdot)\right)$ is continuous in the space $\dot{M}^{2\sigma, p\sigma}(\Rt)$.  We thus have   $\ds{- \frac{1}{-\Delta} \left(  \partial^{\alpha} \text{div}(\V \otimes \U) \right)\in \dot{M}^{2\sigma, p\sigma}(\Rt)}$. 

\medskip

For the second term, by recalling that $\Vert \V \Vert_{L^\infty}=1$,  by using the information $\vec{\nabla}\otimes \V \in \dot{M}^{2,p}(\Rt)\cap L^{\infty}(\Rt)$, and moreover, by Lemma \ref{Lem-U5}, we can write: 
\begin{equation*}
\begin{split}
\Vert \vert \vec{\nabla} \otimes \V \vert^2\, \V \Vert_{\dot{M}^{2\sigma, p\sigma}} \leq  & \,\,  \Vert \vec{\nabla} \otimes \V \Vert_{\dot{M}^{2\sigma, p\sigma}}\Vert \vec{\nabla} \otimes \V \Vert_{L^{\infty}}\Vert \V \Vert_{L^{\infty}} \\
\leq & \,\, c\, \Vert \vec{\nabla} \otimes \V \Vert^{\frac{1}{\sigma}}_{\dot{M}^{2,p}} \Vert \vec{\nabla} \otimes  \V \Vert^{1-\frac{1}{\sigma}}_{L^{\infty}}\, \Vert \vec{\nabla} \otimes \V \Vert_{L^{\infty}}\,\Vert \V \Vert_{L^{\infty}}\\
\leq & \,\, c\, \Vert \vec{\nabla} \otimes \V \Vert^{\frac{1}{\sigma}}_{\dot{M}^{2,p}} \Vert \vec{\nabla} \otimes  \V \Vert^{2-\frac{1}{\sigma}}_{L^{\infty}} < +\infty. 
\end{split}
\end{equation*}
Consequently, by Lemma \ref{Lem-U3} we have  $\ds{\frac{1}{-\Delta} \left( \partial^{\alpha} \left( \vert  \vec{\nabla} \otimes \V  \vert^2\, \V  \right)\right)\in \dot{M}^{2\sigma, p\sigma}(\Rt)}$. 

\medskip

For the third term, always  by  Lemma \ref{Lem-U3} and by  (\ref{Cond-Th-1}), we obtain
\begin{equation}\label{Eq-Force-Ext} 
\left\Vert \frac{1}{-\Delta} \Big( \P\left(  \partial^{\alpha} \text{div}(\mathbb{G})\right) \Big)  \right\Vert_{\dot{M}^{2\sigma, p\sigma}} \leq c \, \left\Vert  \text{div}(\mathbb{G}) \right\Vert_{\dot{M}^{2\sigma, p\sigma}} \leq c \, \Vert \text{div}(\mathbb{G}) \Vert^{\sigma}_{\dot{M}^{2,p}} \, \Vert  \text{div}(\mathbb{G}) \Vert^{1-\sigma}_{L^\infty}<+\infty. 
\end{equation}

We study now the function $\partial^{\alpha} \U $ given in the expression (\ref{eq-der-alpha-U}). For the first term on the right-hand side    we write 
\[ \frac{1}{-\Delta} \left( \P ( \partial^{\alpha} \text{div} (\U \otimes \U) ) \right)= \frac{1}{-\Delta} \left( \P ( \partial^{\alpha_1}  \text{div} \,\partial^{\alpha_2} (\U \otimes \U) ) \right), \quad  \text{where}\quad \vert \alpha_1 \vert=1 \quad \text{and} \quad \vert \alpha_2 \vert =1.\]
Here, we  must verify  that $\partial^{\alpha_2} (\U \otimes \U)  \in \dot{M}^{2\sigma, p\sigma}(\Rt)$. Indeed,  for $i,j=1,2,3$ we write $\partial^{\alpha_2}(U_i U_j)= (\partial^{\alpha_2} U_i) U_j + U_i (\partial^{\alpha_2}U_j)$. Then, as  we have  $\partial^{\alpha_2}\U \in \dot{M}^{2\sigma,p\sigma}(\Rt)$ and $\U \in L^{\infty}(\Rt)$,  we obtain $\partial^{\alpha_2}(\U \otimes \U) \in \dot{M}^{2\sigma, p\sigma}(\Rt)$.  With this information, and the fact that  by Lemma \ref{Lem-U3} the operator $\frac{1}{-\Delta} \left( \P ( \partial^{\alpha_1} \text{div} ( \cdot) ) \right)$ is continuous in the space $\dot{M}^{2\sigma, p\sigma}(\Rt)$, we finally get $\frac{1}{-\Delta} \left( \P ( \partial^{\alpha} \text{div} (\U \otimes \U) ) \right) \in \dot{M}^{2 \sigma, p \sigma}(\Rt)$.

\medskip

The second term on the right-hand  side of (\ref{eq-der-alpha-U}) follows the same ideas above (with $\vec{\nabla}\otimes \V$ instead of $\U$) and we have $\ds{\frac{1}{-\Delta} \left( \P (\partial^{\alpha} \text{div} ( \vec{\nabla} \otimes  \V \odot \vec{\nabla} \otimes \V)) \right)\in \dot{M}^{2 \sigma, p \sigma}(\Rt)}$. 

\medskip 

The third term on the right-hand side is similarly estimated as in (\ref{Eq-Force-Ext}).   \finpv 
\end{enumerate}	

\begin{Lemme}[The iterative process] For all $1 \leq  m \leq k$ and for all $\vert \alpha \vert \leq m$ we assume $\ds{\partial^{\alpha}\V \in \dot{M}^{2\sigma, p \sigma}(\Rt)}$ and $\ds{\partial^{\alpha}\U \in \dot{M}^{2\sigma, p \sigma}(\Rt)}$, for all  $1\leq \sigma < +\infty$. Then,  it  holds true for all $\vert \alpha \vert=k+2$.  
\end{Lemme}	
\pv   We shall follow the main ideas in the proof of the previous lemma.
\begin{enumerate}
\item[$\bullet$] Let $\vert \alpha \vert = k+1$.  We start by studying  the function $\partial^{\alpha}\V$ given in the identity  (\ref{eq-der-V}). Let us verify that each term on the right-hand side belong to the space $\dot{M}^{2\sigma, p\sigma}(\Rt)$. 

\medskip 

For the first term, we split $\alpha=\alpha_1 + \alpha_2$, with $\vert \alpha_1 \vert =1$ and $\vert \alpha_2 \vert=k$. Then, we write
\[ - \frac{1}{-\Delta} \left(  \partial^{\alpha} \text{div}(\V \otimes \U) \right)= - \frac{1}{-\Delta} \left(  \partial^{\alpha_1} \text{div} \,  \partial^{\alpha_2}(\V \otimes \U) \right). \]
In the last expression, we  verify  that $ \partial^{\alpha_2}(\V \otimes \U)  \in \dot{M}^{2\sigma, p\sigma}(\Rt)$. For $i,j=1,2,3$,  by the  Leibinz rule we  write $\ds{ \partial^{\alpha_2} (V_i U_j)= \sum_{\vert \beta \vert \leq k} c_{\alpha_2,\beta} \, \partial^{\beta}  V_i \,   \partial^{\alpha_2 - \beta} U_j}$,  where $c_{\alpha_2,\beta}>0$ is a constant depending on the multi-indices $\alpha_2$ and $\beta$.  Then, as by the recurrence hypothesis  we have $\ds{\partial^{\alpha_2- \beta} V_i \in \dot{M}^{2\sigma, p \sigma}(\Rt)}$ and $\ds{\partial^{\beta} U_j \in \dot{M}^{2\sigma, p \sigma}(\Rt)}$, and moreover, by applying the H\"older inequalities,   we  get  $\ds{\partial^{\alpha_2} (\V\otimes \U)  \in \dot{M}^{2\sigma, p\sigma}(\Rt) }$.  On the other hand, recall that  as $\vert \alpha_1 \vert=1$ then by Lemma \ref{Lem-U3} the operator $\ds{- \frac{1}{-\Delta} \left(  \partial^{\alpha_1} \text{div}(\cdot) \right)}$ is continuous in the space $\dot{M}^{2\sigma, p\sigma}(\Rt)$. Finally,  for $\vert \alpha \vert=k+1$ we obtain  $\ds{ - \frac{1}{-\Delta} \left(  \partial^{\alpha} \text{div}(\V \otimes \U) \right) \in \dot{M}^{2\sigma, p\sigma}(\Rt)}$. 

\medskip

For the second term,   we split now $\alpha=\alpha_1+\alpha_2$, with $\vert \alpha_1 \vert=2$ and $\vert \alpha_2 \vert= k-1$. Then we write
\[ \frac{1}{-\Delta} \left( \partial^{\alpha} \left( \vert  \vec{\nabla} \otimes \V  \vert^2\, \V  \right)\right)= \frac{1}{-\Delta} \left( \partial^{\alpha_1} \left( \partial^{\alpha_2}  \left( \vert  \vec{\nabla} \otimes \V  \vert^2\, \V  \right)\right)\right). \] 
As before, we must verify that $\ds{\partial^{\alpha_2}  \left( \vert  \vec{\nabla} \otimes \V  \vert^2\, \V  \right) \in \dot{M}^{2\sigma, p\sigma}(\Rt)}$. Always by the Leibinz rule\footnote{For the sake of simplicity, we omit the constants.}, for $i,j,k=1,2,3$ we have the following identity $\ds{\partial^{\alpha_2} \left( (\partial_i V_j)^2 V_k \right)= \sum_{\vert \beta \vert \leq k-1}  \partial^{\beta} \left(  (\partial_i V_j)^2 \right) \,  \partial^{\alpha_2 - \beta} V_k}$.  Moreover, in order to compute the term $\ds{\partial^{\beta} \left(  (\partial_i V_j)^2 \right) }$, we make use again  of the Leibinz rule to write $\ds{\partial^{\beta} \left(  (\partial_i V_j)^2 \right) =  \sum_{\vert \gamma \vert \leq \vert \beta\vert} (\partial^{\gamma} \partial_i V_j) \,  (\partial^{\beta-\gamma} \partial_i V_j)}$.  Thus, by gathering these identities we  obtain 
\[ \partial^{\alpha_2}  \left( \vert  \vec{\nabla} \otimes \V  \vert^2\, \V  \right)= \sum_{\vert \beta \vert \leq k-1} \sum_{\vert \gamma \vert \leq \vert \beta \vert}  \partial^{\gamma} (\partial_i V_j) \, \partial^{\beta-\gamma} (\partial_i V_j)  \,  (  \partial^{\alpha_2 - \beta} V_k). \]
 Once we have this identity, by using the recurrence hypothesis and the H\"older inequalities we get $\ds{\ds{\partial^{\alpha_2}  \left( \vert  \vec{\nabla} \otimes \V  \vert^2\, \V  \right) \in \dot{M}^{2\sigma, p\sigma}(\Rt)}}$. Finally, as we have $\vert \alpha _2 \vert=2$, by Lemma \ref{Lem-U3} the operator $\ds{\frac{1}{-\Delta} \left( \partial^{\alpha_1} \left( \cdot \right)\right)}$ is continuous in $\dot{M}^{2\sigma, p\sigma}(\Rt)$ and for $\vert \alpha \vert=k+1$ we get $\ds{\frac{1}{-\Delta} \left( \partial^{\alpha} \left( \vert  \vec{\nabla} \otimes \V  \vert^2\, \V  \right)\right) \in \dot{M}^{2\sigma, p\sigma}(\Rt)}$. 

\medskip

For the third term,   we split again $\alpha=\alpha_1 +\alpha_2$ with $\vert \alpha_1 \vert =1$ and $\vert \alpha_2\vert=k$. Then, by Lemma \ref{Lem-U3} and by  (\ref{Cond-Th-1}) we have:
\begin{equation}\label{Estim-Force-2}
\begin{split}
\left\Vert \frac{1}{-\Delta} \P \left( \partial^{\alpha}\, \text{div} \left( \mathbb{G} \right) \right)\right\Vert_{\dot{M}^{2\sigma, p\sigma}} = &\,  \left\Vert \frac{1}{-\Delta} \P \left( \partial^{\alpha_1}\, \text{div} \, \partial^{\alpha_2} \left( \mathbb{G} \right)\right) \right\Vert_{\dot{M}^{2\sigma, p\sigma}} \leq c\, \left\Vert \partial^{\alpha_2} \mathbb{G} \right\Vert_{\dot{M}^{2\sigma, p\sigma}} \\
\leq & \,  c\, \Vert \partial^{\alpha_2} \mathbb{G} \Vert^{\frac{1}{\sigma}}_{\dot{M}^{2,p}}\, \Vert \partial^{\alpha_2} \mathbb{G} \Vert^{1-\sigma}_{L^{\infty}}<+\infty. 
\end{split}
\end{equation}
Once we have $\partial^{\alpha}\V \in \dot{M}^{2\sigma, p\sigma}(\Rt)$ for all $\vert \alpha \vert \leq k+1$,  by using the identity (\ref{eq-der-alpha-U}) and by following the same estimates above, we obtain  $\partial^\alpha \U \in \dot{M}^{2\sigma, p\sigma}(\Rt)$ for $\vert \alpha \vert=k+1$.

\item[$\bullet$] Let $\vert \alpha \vert=k+2$. Once we have $\partial^{\alpha}\V \in \dot{M}^{2\sigma, p\sigma}(\Rt)$ and $\partial^\alpha \U \in \dot{M}^{2\sigma, p\sigma}(\Rt)$ for all $\vert \alpha \vert \leq k+1$,   we just repeat again the estimates above  to obtain  $\partial^{\alpha} \U, \, \partial^{\alpha} \V \in \dot{M}^{2\sigma, p\sigma}(\Rt)$ for $\vert \alpha \vert=k+2$.  \finpv 
\end{enumerate}	
\begin{Remarque}\label{Rmk-tech}  By (\ref{Cond-Th-1}) and by the estimate (\ref{Estim-Force-2}) (a similar estimate holds for $\mathbb{F}$), $k+2$ is  the maximum gain of regularity expected for $\partial^{\alpha} \U$ and $\partial^{\alpha} \V$ given in (\ref{eq-der-alpha-U}) and (\ref{eq-der-V}) respectively.
\end{Remarque}	

It remains to study the pressure term in the first equation of the system  (\ref{EickLes}). For this, we apply the divergence operator in this equation  to obtain that  $P$ is   related to $\U$ and $\vec{\nabla} \otimes \V$ through the expression: 
\begin{equation}\label{Expression-P} 
P =   \frac{1}{-\Delta} \left(  \text{div}\left( \text{div}(\U \otimes \U) \right) \right)+  \frac{1}{-\Delta} \left(  \text{div}\left( \text{div}(\vec{\nabla} \otimes \V \odot \vec{\nabla} \otimes \V  ) \right) \right) + \frac{1}{-\Delta}\left( \text{div} \left( \text{div} (\mathbb{F}) \right) \right).
\end{equation}	 
By following the same estimates in the proof of the previous lemmas, we obtain:
\begin{Lemme} For all $\vert \alpha \vert \leq k+1$  we have $ \partial^{\alpha} P \in \dot{M}^{2\sigma, p\sigma}(\Rt)$.   
\end{Lemme}	
\pv  We study each term  on the right-hand side  in the identity (\ref{Expression-P}). For the first term,  as we have $\partial^{\alpha} \U \in \dot{M}^{2\sigma, p\sigma}$ for all $\vert \alpha \vert \leq k+2$, and as we have $\U \in L^{\infty}(\Rt)$,  then we get $\ds{\frac{1}{-\Delta} \left(  \text{div}\left( \text{div}(U_i\, U_j) \right) \right) \in \dot{M}^{2\sigma, p\sigma}(\Rt)}$ for $\vert \alpha \vert \leq k+1$.  Similarly, for the second term,  as  $\partial^{\alpha} \V \in \dot{M}^{2\sigma, p\sigma}$ for all $\vert \alpha \vert \leq k+2$, and as $\vec{\nabla} \otimes  \in L^{\infty}(\Rt)$,  we get $\ds{\frac{1}{-\Delta} \left(  \text{div}\left( \text{div}(\vec{\nabla} \otimes \V \odot \vec{\nabla} \otimes \V  ) \right) \right) \in \dot{M}^{2\sigma, p\sigma}(\Rt)}$ for $\vert \alpha \vert \leq k+1$. Finally,  by (\ref{Cond-Th-1}) for the third term  we have $\ds{\frac{1}{-\Delta}\left( \text{div} \left( \text{div} (\mathbb{F}) \right) \right) \in \dot{M}^{2\sigma, p\sigma}(\Rt)}$, for all $\vert \alpha \vert \leq k+1$.  \finpv 

\medskip

We are able to finish the proof of Theorem \ref{Th}. Recall that we have the embedding  $\ds{\dot{M}^{2,p}(\Rt) \subset \dot{M}^{1,p}(\Rt)}$. Then, for all   $\vert \alpha  \vert \leq k+1$, by Lemma \ref{Lem-U4} the functions $\partial^{\alpha} \U$  and $\partial^{\alpha} \V$  are H\"older continuous with exponent $0<1-3/p<1$. Moreover,  for $\vert \alpha \vert \leq k$ the function $\partial^{\alpha} P$  is also  H\"older continuous with the same exponent.   Theorem \ref{Th} is now proven. \finpv

\section{Proof of Theorem \ref{Th2}}  

The proof of this theorem  essentially follows the same lines in the proof of  Theorem \ref{Th}. Consequently,  we will give a sketch of the main steps.  We consider first the evolution problem for the (MHD) equations: 
\begin{equation}\label{MHD-Ev}
\begin{cases}\vspace{2mm}
\partial_t \vu - \Delta \vu + \P ( \text{div}(\vu \otimes \vu)) - \P ( \text{div}( \b \otimes \b ) ) =\P\, (\text{div}(\mathbb{F})),   \\ \vspace{2mm}
\partial_t \b - \Delta \b + \text{div}(\b \otimes \vu) - \text{div}(\vu \otimes \b)=\text{div}(\mathbb{G}),  \\ \vspace{2mm}
\text{div}(\vu)=\text{div}(\b)=0,\\
\vu(0,\cdot)= \vu_0\, \quad \b(0,\cdot)=\b_0. 
\end{cases}
\end{equation} 
By following the  same computations done in the proof of  Proposition \ref{Prop1},  for the initial data $(\vu_0, \b_0)\in \dot{M}^{2,p}(\Rt)$ (with $p>3$) there exists a time $T_0>0$ and  there exists $(\vu, \b ) \in  \mathcal{C}_{*}([0,T_0], \dot{M}^{2,p}(\Rt))$ a unique solution of (\ref{MHD-Ev}) which verifies: 
\begin{equation}\label{EqLinf}
\sup_{0 < t <T} t^{\frac{3}{2p}} \Big( \Vert \vu(t,\cdot)\Vert_{L^\infty} + \Vert \b(t,\cdot) \Vert_{L^\infty}\Big)<+\infty.  
\end{equation}

As above, the key idea is the fact that the time-independent functions $(\U, \B) \in \dot{M}^{2,p}(\Rt)$, which solve the system (\ref{MHD}), also belong to the space  $\mathcal{C}_{*}([0,T], \dot{M}^{2,p}(\Rt))$ and they solve the evolutionary system (\ref{MHD-Ev}) with the initial data $\vu_0 = \U$ and $\b_0=\B$. Then, by uniqueness of solutions in this space we obtain that $\U$ and $\B$ verify (\ref{EqLinf}), hence, we have that $\U \in L^{\infty}(\Rt)$ and $\B \in L^{\infty}(\Rt)$. 

\medskip

On the other hand,  $(\U, \B) \in \dot{M}^{2,p}(\Rt)$ verify the following integral equations:
\begin{equation}\label{MHD-Int}
\begin{array}{ll}\vspace{2mm} 
\ds{\U = - \frac{1}{-\Delta} \left( \P ( \text{div} (\U \otimes \U) ) \right) - \frac{1}{-\Delta} \left( \P ( \text{div} ( \B \otimes  \B )) \right)} + \frac{1}{-\Delta}\left( \P(\text{div}(\mathbb{F}))\right), \\ \vspace{2mm} 
\ds{\B =- \frac{1}{-\Delta} \left( \text{div}(\B \otimes \U) \right) - \frac{1}{-\Delta} \left(\text{div}(\U \otimes \B ) \right)}+\frac{1}{-\Delta}\left( \text{div}(\mathbb{G})\right), \\ 
\vspace{2mm} 
\ds{ P= \sum_{i,j=1}^{3} \mathcal{R}_i \mathcal{R}_j \Big( U_i U_j + B_i B_j+ F_{i,j} \Big) }, \quad \mathcal{R}_i = \frac{\partial_i}{\sqrt{-\Delta}}, \,\, i=1,2,3.
\end{array} 
\end{equation}
As we have  $\U \in L^{\infty}(\Rt)$ and $\B \in L^{\infty}(\Rt)$, and moreover, by  following the computations performed  in the proof of Proposition \ref{Prop3},  for all  $\vert \alpha \vert \leq k+2$ we have  $\partial^{\alpha}\U \in \dot{M}^{2,p}(\Rt)$, $\partial^\alpha \B \in \dot{M}^{2,p}(\Rt)$;  and for all $\vert \alpha \vert \leq k+1$ we have  $\partial^{\alpha} P \in \dot{M}^{2,p}(\Rt)$.  Theorem \ref{Th2} is proven.  \finpv

\section{Estimates in the Gevrey class}
\subsection{Proof of  Theorem \ref{Th3}}  
We consider the evolution problem for the (MHD) equations given in (\ref{MHD-Ev}),  but with  external forces $\f$ and $\g$ in the first and the second equation respectively:  
\begin{equation}\label{MHD-Ev-2}
\begin{cases}\vspace{2mm}
\partial_t \vu - \Delta \vu + \P ( \text{div}(\vu \otimes \vu)) - \P ( \text{div}( \b \otimes \b ) ) = \f,   \\ \vspace{2mm}
\partial_t \b - \Delta \b + \text{div}(\b \otimes \vu) - \text{div}(\vu \otimes \b)=\g,  \\ \vspace{2mm}
\text{div}(\vu)=\text{div}(\b)=0,\\
\vu(0,\cdot)= \vu_0\, \quad \b(0,\cdot)=\b_0. 
\end{cases}
\end{equation} 
Then, we recall the following  classical result: 
\begin{Lemme}\label{Lem1} Let $(\vu_0, \b_0) \in \dot{H}^1(\Rt)$ be a divergence-free initial data. Moreover, for $0<T<+\infty$ we assume  $\f , \g \in \mathcal{C}([0,T], \dot{H}^{1}(\Rt))$. Then, there exists a time $0<T_0<T$, depending on $\vu_0$, $\b_0$ $\f$ and $\g$, and there exists  $(\vu, \b) \in \mathcal{C}([0,T_0], \dot{H}^{1}(\Rt))$ which is the unique  solution of (\ref{MHD-Ev-2}). 
\end{Lemme}	
We refer to Theorem $7.1$ of the book \cite{PGLR1} for a proof in the case of the Navier-Stokes equations (with $\b=0$). However, this proof can be easily adapted for the (MHD) system  (\ref{MHD-Ev-2}). 

\medskip

We recall now that for $b >0$ and for $t>0$,  the weighted exponential operator $e^{b \, \sqrt{t}   \sqrt{-\Delta}}$ is defined in  the Fourier variable as  $\ds{\mathcal{F}\Big( e^{b \, \sqrt{t}  \sqrt{-\Delta}} \varphi(t,\cdot)  \Big)(\xi)= e^{b\, \sqrt{t} \vert \xi \vert} \,  \widehat{\varphi}(t,\xi)}$, for all  $\varphi \in \mathcal{S}([0,+\infty[ \times \Rt)$.  In the  following result we prove that the solution $(\vu, \b)$ obtained above belong to the Gevrey class, provided that 
\begin{equation}\label{Hyp-forces}
e^{b \, \sqrt{t}   \sqrt{-\Delta}} \f \, \in \mathcal{C}([0,T_0], \dot{H}^{1}(\Rt)) \quad \text{and} \quad e^{b \, \sqrt{t}   \sqrt{-\Delta}} \g \in \mathcal{C}([0,T_0], \dot{H}^{1}(\Rt)).
\end{equation}
\begin{Proposition}\label{Prop4}  Assume that $\f$ and $\g$ verify (\ref{Hyp-forces}). Then, there exists a time $0<T_1<T_0$ such that the unique solution $(\vu, \b) \in \mathcal{C}([0,T_0], \dot{H}^{1}(\Rt))$ of (\ref{MHD-Ev-2}) obtained in Lemma \ref{Lem1} verifies:
\begin{equation}\label{COND}
e^{b\, \sqrt{t}\, \sqrt{-\Delta}} \,  \vu \in \mathcal{C}(]0,T_1[, \dot{H}^{1}(\Rt)) \quad \text{and} \quad  e^{b\, \sqrt{t}\, \sqrt{-\Delta}} \, \b \in \mathcal{C}(]0,T_1[, \dot{H}^{1}(\Rt)). 
\end{equation}	
\end{Proposition}	
\pv  We split the proof in two main steps.
\begin{enumerate}
\item[$\bullet$]  For a time $0<T_1<T_0$ (which we shall set small enough)  we consider the  space 
\[ E= \left\{  f \in \mathcal{C}([0, T_1[, \dot{H}^{1}(\Rt)) : \,\, e^{b \, \sqrt{t} \sqrt{-\Delta}} f \in \mathcal{C}(]0, T [, \dot{H}^{1}(\Rt)) \right\}, \]
which is a Banach space with the norm  $\ds{\Vert \cdot \Vert_{E}= \left\Vert  e^{b \, \sqrt{t} \sqrt{-\Delta}} (\cdot) \right\Vert_{L^{\infty}_{t} \dot{H}^{1}_{x}}}$. For the initial data  $(\vu_0, \b_0) \in \dot{H}^{1}(\Rt)$, we will construct a mild  solution $(\vu_1, \b_1)$ of the equation  (\ref{MHD-Ev-2}) in the space $E$.   We recall that this mild solution solves the equations:
\begin{equation}\label{equ2}
\begin{split}
\vu_1(t,\cdot)=&\,\, e^{t \Delta}\vu_0 +\int_{0}^{t} e^{(t-s)\Delta }\,\f \, ds + \underbrace{\int_{0}^{t} e^{(t-s)\Delta }\,\P \, \text{div} \left( \vu_1 \otimes \vu_1 \right)(s,\cdot)\, ds}_{{\bf B_1}(\vu_1,\vu_1)}\\
&\,\, - \underbrace{\int_{0}^{t} e^{(t-s)\Delta }\,\P \, \text{div} \left(\b_1 \otimes \b_1 \right)(s,\cdot)\, ds}_{{\bf B_2}(\b_1,\b_1)},
\end{split}
\end{equation}
\begin{equation}\label{eqb}
\begin{split}
\b_1(t,\cdot)= &\,\, e^{t\Delta} \b_0  +\int_{0}^{t} e^{(t-s)\Delta }\,\g\, ds+ \underbrace{\int_{0}^{t} e^{(t-s)\Delta}\, \text{div}\left(\b_1  \otimes  \vu_1\right)(s,\cdot)\, ds}_{{\bf B_3}(\b_1,\vu_1)}\\
&\,\,- \underbrace{\int_{0}^{t} e^{(t-s)\Delta}\, \text{div}\left(\vu_1 \otimes \b_1\right)(s,\cdot)\, ds}_{{\bf B_4}(\vu_1,\b_1)}. 
\end{split}
\end{equation}
 The linear terms in the system (\ref{equ2})-(\ref{eqb}) are easy to estimate and  for a constant $c>0$, depending on $b$ and $T$, we have: 
\begin{equation}\label{lin}
\begin{split}
&\,\,\left\Vert \left( e^{t\, \Delta} \vu_0 +  \int_{0}^{t} e^{(t-s)\Delta }\,\f \, ds, \,\,  e^{t\, \Delta} \b_0 + \int_{0}^{t} e^{(t-s)\Delta }\,\g \, ds\right) \right\Vert_{E} \\
\leq &\,\,   c  \left( \Vert \vu_0 \Vert_{\dot{H}^{1}}+ \Vert \b_0 \Vert_{\dot{H}^{1}} + \left\Vert e^{b\,\sqrt{t}\sqrt{-\Delta}} \f \right\Vert_{L^{\infty}_{t} \dot{H}^{1}_{x}} + \left\Vert e^{b\,\sqrt{t}\sqrt{-\Delta}} \g \right\Vert_{L^{\infty}_{t} \dot{H}^{1}_{x}}\right).
\end{split}
\end{equation}
On the other hand, as the four nonlinear terms ${\bf B_i}(\cdot,\cdot)$, with  $i=1,\cdots, 4$, in the system  (\ref{equ2})-(\ref{eqb}) have the same structure,  we shall  only estimate the first one.  By the Plancherel's formula  we can write: 
\begin{equation}\label{estim}
\begin{split}
\Vert {\bf B_1}(\vu_1, \vu_1) \Vert_E = & \sup_{0\leq t \leq T} \left\Vert e^{b\,\sqrt{t} \sqrt{-\Delta}} \left(  \int_{0}^{t} e^{(t-s)\Delta} \P (\text{div} (\vu_1 \otimes \vu_1) )(s,\cdot) ds \right)  \right\Vert_{\dot{H}^1}  \\
\leq &  \sup_{0\leq t \leq T} \left\Vert \vert \xi \vert \, e^{b \sqrt{t} \vert \xi \vert}\, \int_{0}^{t} e^{-(t-s)\vert \xi \vert^2} \, \vert \xi \vert \left(  \widehat{\vu_1} \ast \widehat{\vu_1} \right) (s,\cdot) ds \right\Vert_{L^2} \\
\leq &  \sup_{0\leq t \leq T} \left\Vert \vert \xi \vert^2 \,  \int_{0}^{t} e^{-(t-s)\vert \xi \vert^2} \, e^{b \sqrt{t} \vert \xi \vert}\, \left(  \widehat{\vu_1} \ast \widehat{\vu_1} \right) (s,\cdot) ds \right\Vert_{L^2}. 
\end{split}
\end{equation}
We must study now the expression $\ds{e^{b \sqrt{t} \vert \xi \vert} \left(  \widehat{\vu} \ast \widehat{\vu} \right) (s,\xi)}$.  We remark that for all $\xi, \eta \in \Rt$  we have $\ds{e^{b \sqrt{t} \vert \xi \vert} \leq e^{b \sqrt{t} \vert \xi - \eta \vert} \, e^{b \sqrt{t} \vert \eta \vert}}$. Then, we obtain the following pointwise estimate: 
\[ \left\vert  e^{b \sqrt{t} \vert \xi \vert} \left(  \widehat{\vu_1} \ast \widehat{\vu_1} \right) (s,\xi) \right\vert \leq \left(  \left(  e^{b \sqrt{t} \vert \xi \vert} \vert \widehat{\vu_1} \vert  \right) \ast  \left(  e^{b \sqrt{t} \vert \xi \vert} \vert \widehat{\vu_1} \vert  \right) \right)(s,\xi).\] 
By getting back to the estimate (\ref{estim})  we write
\begin{equation}
\begin{split}
& \sup_{0\leq t \leq T} \left\Vert \vert \xi \vert^2 \, \int_{0}^{t} e^{-(t-s)\vert \xi \vert^2} e^{b \sqrt{t} \vert \xi \vert} \left(  \widehat{\vu_1} \ast \widehat{\vu_1} \right) (s,\cdot) ds \right\Vert_{L^2} \\
\leq & \sup_{0\leq t \leq T} \left\Vert \vert \xi \vert^2 \, \int_{0}^{t} e^{-(t-s)\vert \xi \vert^2} \left(  \left(  e^{b \sqrt{t}  \vert \xi \vert} \vert \widehat{\vu_1} \vert  \right) \ast  \left(  e^{b \sqrt{t}  \vert \xi \vert} \vert \widehat{\vu_1} \vert  \right) \right)(s,\cdot) ds \right\Vert_{L^2} \\
\leq & \sup_{0\leq t \leq T} \int_{0}^{t} \left\Vert \vert \xi \vert^{3/2} \, e^{-(t-s)\vert \xi \vert^2} \, \vert \xi \vert^{1/2}\left(  \left(  e^{b \sqrt{t}  \vert \xi \vert} \vert \widehat{\vu_1} \vert  \right) \ast  \left(  e^{b \sqrt{t}  \vert \xi \vert} \vert \widehat{\vu_1} \vert  \right) \right)(s,\cdot) ds \right\Vert_{L^2}=(a). 
\end{split}
\end{equation}  
By applying again  the Plancherel's formula we get back to the spatial variable. Moreover, by the well-known properties of the heat kernel and by the product laws in the homogeneous Sobolev spaces, we have:
\begin{equation*}
\begin{split}
(a) \leq & \,\, \sup_{0\leq t \leq T}  \, \int_{0}^{t} \left\Vert (-\Delta)^{3/4} h_{(t-s)}(\cdot) \ast  (-\Delta)^{1/2} \left(  \left(  e^{b \sqrt{t} \vert \xi \vert} \vert \widehat{\vu_1} \vert  \right)^{\vee} \times  \left(  e^{b \sqrt{t}  \vert \xi \vert} \vert \widehat{\vu_1} \vert  \right)^{\vee} \right)(s,\cdot) \right\Vert_{L^2}\, ds  \\
\leq & \,\, \sup_{0\leq t \leq T}  \, \int_{0}^{t} \left\Vert (-\Delta)^{3/4} h_{(t-s)}(\cdot)  \right\Vert_{L^1} \, \left\Vert  (-\Delta)^{1/2} \left(  \left(  e^{b \sqrt{t} \vert \xi \vert} \vert \widehat{\vu_1} \vert  \right)^{\vee} \times \left(  e^{b \sqrt{t}  \vert \xi \vert} \vert \widehat{\vu_1} \vert  \right)^{\vee} \right)(s,\cdot)  \right\Vert_{L^2} \, ds \\
\leq &\,\, c \,  T^{1/4} \, \sup_{0\leq s \leq T}  \,  \left\Vert \left(  \left(  e^{b \sqrt{t} \vert \xi \vert} \vert \widehat{\vu_1} \vert  \right)^{\vee} \times  \left(  e^{b \sqrt{t}  \vert \xi \vert} \vert \widehat{\vu_1} \vert  \right)^{\vee} \right)(s,\cdot) \right\Vert_{\dot{H}^{1/2}} \\
\leq &\,\, c \,  T^{1/4} \, \sup_{0\leq s \leq T}  \,  \left\Vert \left(  \left(  e^{b \sqrt{t} \vert \xi \vert} \vert \widehat{\vu_1} \vert  \right)^{\vee}  \right\Vert_{\dot{H}^{1}} \, \left\Vert  \left(  e^{b \sqrt{t}  \vert \xi \vert} \vert \widehat{\vu_1} \vert  \right)^{\vee} \right)(s,\cdot) \right\Vert_{\dot{H}^{1}}\\
\leq &\,\,  c\, T^{1/4}\, \Vert \vu_1 \Vert_{E}\, \Vert \vu_1 \Vert_{E}. 
\end{split}
\end{equation*}
The other nonlinear terms are treated in the same way; and we can write: 
\begin{equation}\label{non-lin}
\Vert {\bf B_1 }(\vu_1, \vu_1)  \Vert_{E}+ \Vert {\bf B_2} (\b_1, \b_1)  \Vert_{E}+ \Vert {\bf B_3 }(\b_1, \vu_1)  \Vert_{E}+ \Vert {\bf B_4 }(\vu_1, \b_1) \Vert_{E} \leq\,\, c\, T^{1/4}\, \Vert (\vu_1, \b_1) \Vert^{2}_{E}. 
\end{equation}
Thus, by (\ref{lin}) and (\ref{non-lin}), we fix a time $0<T_1<T_0$ small enough such that we cal apply the Picard's iterative schema to obtain a solution $(\vu_1, \b_1) \in E$ of the system (\ref{equ2})-(\ref{eqb}). 

\item[$\bullet$]  By recalling that $E \subset \mathcal{C}([0,T_0], \dot{H}^{1}(\Rt))$, and moreover, by uniqueness of solutions in the space $\mathcal{C}([0,T_0], \dot{H}^{1}(\Rt))$, for  all $0<t<T_1$ we have the identity $\left(\vu_1(t,\cdot) , \b_1(t,\cdot)\right)= \left(\vu(t,\cdot), \b(t,\cdot)\right)$. Consequently, the solution $(\vu, \b)$ obtained in Lemma \ref{Lem1} verifies (\ref{COND}). Proposition \ref{Prop4} is proven. \finpv
\end{enumerate}	

\medskip 

We are able to finish the proof of  Theorem \ref{Th3}. We observe that the time-independent solution $(\U, \B) \in \dot{H}^{1}(\Rt)$ of the system (\ref{MHD}) verifies $(\U, \B) \in \mathcal{C}([0,T_1], \dot{H}^{1}(\Rt))$. Moreover, this solution also solves the evolutionary problem (\ref{MHD-Ev-2}) with initial data $(\vu_0, \b_0)= (\U, \B)$ and with the external forces  $\f=\text{div}(\mathbb{F})$ and $\g=\text{div}(\mathbb{G})$.  

\medskip 

As $\mathbb{F}, \mathbb{G} \in G^{0}_{b}(\Rt)$ with $b>0$, we define the parameter  $\beta=  \frac{2  b}{3 \sqrt{T_0}}>0$,  and then we rewrite  the external forces $\f$ and $\g$ as follows:
\begin{equation}\label{Forces}
\f = e^{-\beta\sqrt{t} \sqrt{-\Delta}} \left( e^{\beta\sqrt{t} \sqrt{-\Delta}}\, \text{div}(\mathbb{F}) \right), \quad  \g = e^{-\beta\sqrt{t} \sqrt{-\Delta}} \left( e^{\beta\sqrt{t} \sqrt{-\Delta}}\, \text{div}(\mathbb{G}) \right).
\end{equation}
Then, we shall prove that $\f, \g \in \mathcal{C}([0,T_0], \dot{H}^{1}(\Rt))$, and moreover,  we shall prove that $\f$ and $\g$ verify (\ref{Hyp-forces}). Indeed, we write:
\begin{equation*}
\begin{split}
&\,\,\sup_{0<t<T_0} \left\Vert \left( e^{\beta\sqrt{t}\sqrt{-\Delta}} \f, \,\ e^{\beta\sqrt{t}\sqrt{-\Delta}} \g\right) \right\Vert^{2}_{\dot{H}^{1}}\\
 \leq  &\,\,  \sup_{0<t<T_0} \, \int_{\Rt} \vert \xi \vert^2\, e^{2 \beta \sqrt{t} \vert \xi \vert} \left( \vert \widehat{f}(t,\xi)\vert^2 + \vert \widehat{g}(t,\xi)\vert^2\right)\, d\xi \\
\leq&\,\,  \sup_{0<t<T_0} \, \int_{\Rt} \vert \xi \vert^2\, e^{2 \beta \sqrt{t} \vert \xi \vert}  \left(\vert \widehat{\text{div}(\mathbb{F})}(\xi) \vert^2 + \vert\widehat{\text{div}(\mathbb{G})}(\xi)\vert^2 \right)\, d\xi \\
\leq &\,\,  \sup_{0<t<T_0} \, \int_{\Rt} \vert \xi \vert^4\, e^{2 \beta\sqrt{t} \vert \xi \vert} \left(\vert \widehat{\mathbb{F}}(\xi) \vert^2+\vert \widehat{\mathbb{G}}(\xi) \vert^2\right)\, d\xi \\
\leq & \,\, \frac{1}{(\sqrt{T_0}\, \beta)^4} \int_{\Rt} (\sqrt{T_0} \beta \vert \xi \vert)^4 \, e^{2 (\sqrt{T_0} \beta \vert \xi \vert)} \, \left( \vert \widehat{\mathbb{F}}(\xi) \vert^2+\vert \widehat{\mathbb{G}}(\xi) \vert^2\right)\, d\xi\\
\leq & \,\, \frac{1}{(\sqrt{T_0}\,\beta)^4} \int_{\Rt} \, e^{3 \sqrt{T_0} \beta \vert \xi \vert} \, \left(\vert \widehat{\mathbb{F}}(\xi) \vert^2+\vert \widehat{\mathbb{G}}(\xi) \vert^2\right)\, d\xi \\
\leq & \,\, c\, \int_{\Rt} \, e^{2 b \vert \xi \vert} \, \left(\vert \widehat{\mathbb{F}}(\xi) \vert^2+\vert \widehat{\mathbb{G}}(\xi) \vert^2\right)\, d\xi  <+\infty. 
\end{split}
\end{equation*}

\medskip 

By Lemma \ref{Lem1} we have that $(\U, \B)$ is the unique solution of the system (\ref{MHD-Ev-2}) in the space $\mathcal{C}([0,T_0], \dot{H}^{1}(\Rt))$. Thereafter, by Proposition \ref{Prop4} we have $\ds{e^{\beta\, t \sqrt{-\Delta}} \U \in \mathcal{C}(]0,T_1[, \dot{H}^{1}(\Rt))}$ and $\ds{ e^{\beta\, t \sqrt{-\Delta}} \B \in \mathcal{C}(]0,T_1[, \dot{H}^{1}(\Rt))}$. We thus set $b_1= \beta\frac{T_1}{2}>0$ and we get $\U \in G^{1}_{b_1}(\Rt)$ and $\B \in G^{1}_{b_1}(\Rt)$. Finally, by the third identity in (\ref{MHD-Int})  and by the product laws in the homogeneous Sobolev spaces we have $P \in G^{1/2}_{b_1}(\Rt)+G^{0}_{b}(\Rt)$. Theorem \ref{Th3} is proven. \finpv

\subsection{Proof of Corollary \ref{Col3}} 
In the case $\mathbb{F}=\mathbb{G}=0$, the external forces $\f$ and $\g$ given in (\ref{Forces}) are null and evidently  they verify (\ref{Hyp-forces}) for any $b>0$. Then, the subsequent parameters  $\beta=  \frac{2  b}{3 \sqrt{T_0}}>0$ and $b_1= \beta \frac{T_1}{2}>0$ can be fixed arbitrary. \finpv

\section*{Data availability statement}
Data sharing not applicable to this article as no datasets were generated or analyzed during the current study.


\begin{thebibliography}{40}
\bibitem{DCh} D. Chamorro. \emph{Espacios de Lebesgue y de Lorentz.} Vol. 3. hal-01801025v1 (2018). 


\bibitem{Ericksen} J.L. Ericksen. \emph{Hydrostatic theory of liquid crystals}. Arch. Rational Mech. Anal, 9:371–378 (1962) 
\bibitem{PFOJ1} P.G. Fernández-Dalgo \& O. Jarrín. \emph{Existence of infinite-energy and discretely self-similar global weak solutions for 3D MHD equations}. arXiv:1910.11267. To appear in the Journal of Mathematical Fluid Mechanics (2020).

\bibitem{Foias} C. Foias \&  R. Temam. \emph{Gevrey class regularity for the solutions of the Navier-Stokes equations}.
Journal of Functional Analysis, 87(2):359–369 (1989). 

\bibitem{Huang} T. Huang. \emph{Regularity and uniqueness for a class of solutions to the hydrodynamic flow of nematic liquid crystals}  Analysis and Applications Vol. 14, No. 04, pp. 523-536 (2016). 

\bibitem{Huang2} T. Huang \& C. Y. Wang. \emph{On uniqueness of heat flow of harmonic maps}. Preprint, arXiv: 1208.1470.


\bibitem{Gennes} P.G. de Gennes. \emph{The physics of liquid crystals}. Oxford University Press, Oxford (1974). 
\bibitem{GigaMiyakawa} Y. Giga \& T. Miyakawa. \emph{Navier-stokes flow in $\Rt$ with measures as initial vorticity and morrey spaces}. Communications in Partial Differential Equations, 14:5, 577-618 (1989).
	\bibitem{HaLiZ} Y. Hao, X. Liu \& X. Zhang. \emph{Liouville theorem for steady-state solutions of simplified Ericksen-Leslie system.} arXiv:1906.06318v1 (2019).  
\bibitem{OJ} O. Jarr\'in. \emph{Liouville theorems for a stationary and non-stationary coupled system of liquid crystal flows}.  Journal of Mathematical Fluid Mechanics, Volume 24, Article number: 50 (2022). 

\bibitem{Kato}T. Kato. Strong Solutions of the Navier-Stokes Equation in Morrey Spaces. Bol. Soc. Bras. Mat., Vol. 22, No. 2: 127-155 (1992).
 \bibitem{KochTataru} H. Koch \& D. Tataru. \emph{Well-posedness for the Navier-Stokes equations} .Adv. Math., 157(1):22–35 (2001).
	\bibitem{LinWang} F.H. Lin \& C.Y. Wang. \emph{Global existence of weak solutions of the nematic liquid crystal flow in dimension three}.  Comm. Pure Appl. Math. 69(8), 1532–1571 (2016). 
\bibitem{PGLR} P.G. Lemari\'e-Rieusset. \emph{Recent developments in the Navier-Stokes problem}, Chapman \& Hall/CRC, (2002).

\bibitem{PGLR1} P.G. Lemari\'e-Rieusset. \emph{The Navier-Stokes Problem in the 21st Century}, Chapman \& Hall/CRC, (2016).
\bibitem{Lin} F.H. Lin. \emph{Nonlinear theory of defects in nematic liquid crystals; phase transition and flow phenomena} Comm. Pure Appl. Math. 42(6): 789–814 (1989). 
\bibitem{LinLiu} F. H. Lin \&  C. Liu. \emph{Partial regularity of the dynamic system modeling the flow of liquid crystals}.  Dis. Cont. Dyn. Syst. 2(1),1-22 (1998). 
\bibitem{LinLinWang}  F. H. Lin, J. Y. Lin \& C. Y. Wang. \emph{Liquid crystal flows in two dimensions}. Arch. Ration. Mech. Anal. 197, no. 1, 297-336 (2010).
	\bibitem{Les} F.M. Leslie. \emph{Some constitutive equations for liquid crystals}. Archive for Rational Mechanics and Analysis,  28(4):265–283 (1968). 

\bibitem{Schnack}  D. D. Schnack. \emph{Lectures in magnetohydrodynamics. With an appendix on extended MHD.} Lecture Notes in Physics, 780. Springer-Verlag, Berlin, (2009).

\bibitem{Wang} C. Y. Wang. \emph{Well-posedness for the heat flow of harmonic maps and the liquid crystal flow with rough initial data}. Arch. Ration. Mech. Anal. 200 no. 1, 1-19 (2011). 

\bibitem{Wang2} C. Y. Wang. \emph{Heat flow of harmonic maps whose gradients belong to $L^{n}_{x}L^{\infty}_{t}$}  Arch. Rational
Mech. Anal., 188:  309-349 (2008). 
\end{thebibliography}
\end{document}